\documentstyle[12pt]{article}
\def\C{{\rm C  \! \! \! l \ }}
\def\R{{\rm I   \! R \ }}

\begin{document}
\begin{flushright} 
{ math.QA/0209158}
\end{flushright}

\begin{center}
{\bf \large ON CONTRACTIONS OF QUANTUM ORTHOGONAL GROUPS} \\
{N.A.Gromov},
{I.V. Kostyakov},  %who made a talk should be underlined
{V.V. Kuratov}                                  \\
{Department of Mathematics, Syktyvkar Branch of IMM,
 Kommunisticheskaya st. 24, Syktyvkar, 167000, Russia\\
E-mail: gromov@dm.komisc.ru}
\end{center}

\begin{abstract}

Instead of zero tending parameters of Wigner--In\"on\"u
our approach to a group contractions is based on use of
the nilpotent commutative generators of Pimenov algebra
$ {\bf D}(\iota), $ which is a subalgebra
of  even part of Grassnamm algebra.
The standard Faddeev quantization of the simple groups
is modified in such a way that the quantum analogs of the
nonsemisimple groups are obtained by contractions.

The contracted quantum  groups
are regarded as the algebras of noncommutative functions
generated by elements $ J_{ik}t_{ik},$ where $J_{ik} $
are some products of   generators
of the algebra $ {\bf D}(\iota)$ and $t_{ik}$ are the
noncommutative generators of guantum group.
Possible contractions of quantum orthogonal groups
are regarded in detail. They essentially depend
 on the choice of primitive elements of the
Hopf algebra. All such choices are considered for quantum group
$SO_{q}(N;C)$ and all allowed contractions in Cayley--Klein
scheme are described.

The quantum deformations of the complex kinematical
groups have been investigated as a contractions of
$ SO_q(5;C) $   and have shown that the result is connected
with the behavior of deformation parameter under contraction.
If deformation parameter $q$ remain unchanged, then
the quantum Euclead $ E_q(4;C) $ and Newton  $ N_q(4;C) $ groups
are obtained. If the deformation parameter is transformed,
then one more nonisomorphic quantum deformation of Newton group $N_v(4;C)$
is obtained.
 But there is no quantum analog of the (complex) Galilei group
$G(1,3)$ in both cases.

According to correspondence principle a new physical theory must include
an old one as a particular case. For space-time symmetries
 this principle is
realized as the chain of contractions of  the kinematical groups:
$$
  S^{\pm}(1,3)\stackrel{ K \rightarrow 0}{\longrightarrow}
  P(1,3)\stackrel{c \rightarrow \infty }{\longrightarrow}G(1,3).
$$
As it was mentioned above there is no quantum deformation of the complex
Galilei group in the standard Cayley--Klein scheme, 
therefore it is not possible to construct the quantum
analog of the full chain of contractions of the (1+3) kinematical groups  
even at the level of a complex groups.

\end{abstract}

\section{Introduction}

Contraction of Lie groups (algebras) is the method of obtaining a new Lie
groups (algebras) from some initial one's with the help of passage to the
limit \cite{IW}. One may define contraction of algebraic structure $(M,*)$ as a
map $\phi_{\epsilon} : (M,*) \rightarrow (N,*')$,
where $(N,*')$ is the algebraic structure of the same type,
isomorphic to $(M,*)$ for  $ \epsilon \neq 0 $ and nonisomorphic to
the initial algebraic structure for $ \epsilon = 0. $
Except for Lie group (algebra) contractions,  graded
contractions \cite{P}, \cite{MP}
 are known, which preserve the grading of Lie algebra.
 Under contractions of bialgebra \cite{VG}
 Lie algebra structure and cocommutator are conserved.
 Hopf algebra (or quantum group \cite{Fad-89}) contractions
are introduced (on the level of quantum algebra \cite{C1}, \cite{C2}
and on the level of quantum group \cite{G})
in such a way that in the limit $ \epsilon \rightarrow 0 $
a new expressions for coproduct, counit and antipode
are consistent with the Hopf algebra axioms.
Recently contractions of the algebraic structures with
bilinear products of arbitrary nature on sections of
finite-dimensional vector bundles was presented \cite{CGM}
and contractions of Lie algebroids and Poisson brackets
was given as an example.

Low dimensional quantum groups have been studed in details.
The two-dimensional Euclidean quantum group $E_q(2)$
was obtained by contractions of the unitary quantum
group $SU_q(2)$ with untouched deformation parameter $q$
in \cite{W}, \cite{SWZ}, \cite{G} and by contractions of the
orthogonal quantum group $ SO_q(3) $ with transformed
deformation parameter in \cite{BHOS}--\cite{Z}.
The quantum Heisenberg group $H_q(1)$ was regarded
as a contraction of $SU_q(2)$ in \cite{ES} and of
$SO_q(3)$ in \cite{BHOS}. A contraction procedure
starting from $SO_q(4)$ was used in \cite{C2} to
determine $E_q(3).$ A contraction of the de Sitter
quantum group leading to a Poincare quantum group
in any dimensions was proposed in \cite{Z1}.
Quantum deformations of the inhomogeneous Lie groups
have been studied in any dimensions by using the projective
(not contraction) method of \cite{SWW}, \cite{CL} for
the multiparametric quantum groups as well \cite{AC}, \cite{AC1}.
%But we did not find in literature a systematic investigations
%of the contractions of quantum groups (not quantum algebras)
%in arbitrary dimensions.
On the other hand $SO_q(3)$ and $SO_q(4)$  are not typical
representatives of the quantum orthogonal groups $SO_q(N)$
for $N=2n+1$ and $N=2n,$ respectively. We shall see that the
number of the allowed contractions for $SO_q(N)$ with  the
transformed deformation parameter is less then the whole
number of the contraction parameters. $SO_q(3)$ and $SO_q(4)$
quantum groups are an exceptions, because both such numbers
are equal (two and three, respectively). Therefore, the
investigation of the  contractions of the quantum
orthogonal groups $SO_q(N)$ for an arbitrary $N$ seems
to worth attention. In present paper  contractions
of the standard deformed quantum  group $ SO_q(N) $ \cite{Fad-89}
are studed in the Cayley--Klein scheme.
The preliminary results for the particular case of identical
permutation was published in \cite{GKK}.

Contractions as a passage to  limit are corresponded with a
physical intuition. At the same time it is desirable
to investigate contractions of an algebraic structures
with the help of pure algebraic tools.
Sometimes it facilitate an investigations,
especially in complicate cases.
It is possible for
classical and quantum Lie groups and algebras if one take into
consideration an algebra $ {\bf D}(\iota) $ with nilpotent
commutative generators. In particular, a motion groups of a
constant curvature spaces (or Cayley-Klein groups) may be obtained
from a classical orthogonal group by replacement its matrix elements
with the  specific  elements  of  the  algebra  $  {\bf  D}(\iota)$
\cite{2}. It is worth to note, that at any stage one may to come back
to the standart In{\"o}n{\"u}--Wigner contraction by putting
an appropriate parameter  tends to zero instead of takes nilpotent value.
In present paper
  the groups under consideration are regarded
according to \cite{Fad-89} as an algebra of noncommutative functions,
but with nilpotent generators.
%From the contraction viewpoint
%Hopf algebra structure of quantum orthogonal group is more rigid
%as compared with the group one.
Possible contractions are essentially
depended on the choice of primitive elements of  Hopf algebra.
We have regarded all variants of such choise for the quantum orthogonal
group $ SO_q(N)$ and for each variant have found all admissible contractions
in Cayley-Klein scheme.

The paper is organized as follows. In Sec. 2, we briefly recall
the matrix realizations of the non-quantum orthogonal Cayley-Klein groups
both in Cartesian and symplectic bases.
 In Sec. 3, the formall definition of
the quantum complex group
$SO_v(N;j;\sigma)$ is given and analysed when the presented structure of
the Hopf algebra is well defined and consistent under nilpotent values
of parameters $j_k$. The results are collected in Theorem 1-4.
The developed approach is applied to the quantum complex kinematic
groups in Sec. 4.
%There is no standart quantum deformation of the Galilei
%$G(1,3)$ group.
The explicit expressions of antipode, coproduct and
relations of $(q,j)$-orthogonality for $SO_v(N;\sigma;j)$ are presented
in Appendices A--C.
We do not pretend to the fullness of the bibliography.
Accessible to us references are included.

\section{Orthogonal Cayley-Klein groups}
\label{sec:level1}

Let us define {\it Pimenov algebra} $ {\bf D}_n(\iota; {\C}) $
as an associative algebra with unit over complex number field
and with nilpotent commutative generators
$ {\iota}_k, \ {\iota}_k^2=0,$
$ {\iota}_k{\iota}_m={\iota}_m{\iota}_k \not =0, \  k \neq m,
\  k,m=1, \ldots, n. $
The general element of ${\bf D}_n(\iota; {\C}) $
is in the form
$$
d=d_0+\sum^{n}_{p=1}\sum_{k_1< \ldots < k_p}d_{k_1\ldots k_p}
{\iota}_{k_1} \ldots {\iota}_{k_p}, \quad  d_0,d_{k_1 \ldots k_p} \in {\C}.
$$
For  $ n=1 $  we have  $ {\bf D}_1({\iota}_1;{\C})  \ni
 d=d_0+d_1{\iota}_1, $ i.e. the elements $d$ are dual (or Study)
numbers when $ d_0,d_1 \in {\bf R} $.
For $ n=2 $ the general element of  $ {\bf D}_2(\iota_1, \iota_2; \C)$ is
$ d=d_0+d_1{\iota}_1+d_2{\iota}_2+d_{12}{\iota}_1{\iota}_2. $
Two elements $ d, \tilde d \in D_n(\iota;{\C})$ are equal if and only if
$ d_0=\tilde d_0, \ d_{k_1 \ldots k_p}=\tilde d_{k_1\ldots k_p},\ p=1,
\ldots,n. $
If $ d=d_k{\iota}_k $ and  $ \tilde d=\tilde d_k{\iota}_k, $ then
the condition $ d=\tilde d, $ which is equivalent to
$ d_k{\iota}_k=\tilde d_k{\iota}_k, $
make possible the consistently definition of the division of nilpotent
generator  $ {\iota}_k $ by itself, namely:
$ {\iota}_k/{\iota}_k=1, \ k=1, \ldots, n. $
Let us stress that the division of different nilpotent generators
$ {\iota}_k/{\iota}_p, \ k \not = p, $ as well as the division of
complex number by nilpotent generators
$ a/{\iota}_k, \ a \in {\C} $  are not defined.
It is convenient to regard the algebras
$ {\bf D}_n(j; {\C}), $ where the parameters
$ j_k=1,\iota_k, \; k=1,\ldots,n. $
If $ m $ parameters are nilpotent $ j_{k_{s}}=\iota_{s}, \; s=1,\ldots,m $
and the other are equal to unit, then we have  Pimenov algebra
$ {\bf D}_m(\iota; {\C}).$

{\it Complex orthogonal Cayley-Klein group} $ SO(N;j;\C) $
is defined as the group of transformations
$ {\xi}'(j)=A(j)\xi(j)$ of complex vector space $O_N(j)$
with Cartesian coordinates
$ \xi^t (j) = (\xi_1,  (1,2)  \xi_2,  \ldots,(1,N)\xi_N)^t, \;$
which preserve the quadratic form
$$
inv(j) =\xi^t(j)\xi(j) = \xi^2_1 + \sum^N_{k=2}({1,k})^2{\xi}^2_k,
$$
where $ j=(j_1, \ldots, j_{N-1}) $, each parameter $ j_k $
takes {\it two} values:
$  j_r=1,{\iota}_r, \ r=1, \ldots, N-1, \ {\xi}_k \in \C$
and
$$
(\mu,\nu)=\prod^{max(\mu,\nu)-1}_{l=min(\mu,\nu)} j_l, \quad (\mu,\mu)=1.
$$
Let us stress, that Cartesian coordinates
of $  O_N(j) $ are special elements of Pimenov algebra
${\bf D}_{N-1}(j; \C). $
It worth notice that the  orthogonal Cayley-Klein groups as well as the unitary
and symplectic   Cayley-Klein   groups   have   been   regarded  in
\cite{G-90}
as the matrix groups with the real matrix elements.
Nevertheless there is a different approach, which gives the same results
for ordinary groups, but is more appropriate from the contraction
quantum group point of view.
According with this approach, the Cayley-Klein group $ SO(N;j;{\C}) $  may be
realised as the matrix group, whose elements are taken from algebra
$ {\bf D}_{N-1}(j;{\C}) $
and in Cartesian basis consist of the $ N \times N $ matrices $A(j)$
with elements
$$
(A(j))_{kp}=(k,p)a_{kp}, \ a_{kp} \in \C.
$$
Matrices $A(j)$ are subject of the additional $j$-orthogonality relations
\begin{equation}
A(j)A^t(j) = A^t(j)A(j)=I.
\label{3}
\end{equation}

Sometimes it is convenient to regard an orthogonal group in so-called
"symplectic" basis. Transformation from Cartesian to symplectic basis
$ x(j)=D\xi (j) $ is made by unitary matrices $D$,
which are a solutions of equation
     \begin{equation}
 D^{t}C_{0}D=I,
     \label{1}
     \end{equation}
where $C_{0}\in M_{N},\; (C_{0})_{ik}= \delta_{ik'},\; k'=N+1-k$
To obtain all solutions of equation (\ref{1}),
 take one of them, namely
$$
D=\frac{1}{\sqrt{2}}
\left ( \begin{array}{cc}
      I &     -i{\tilde C_0} \\
      {\tilde C_0} &   iI
      \end{array} \right ),    \       N=2n,
$$
\begin{equation}
D=\frac{1}{\sqrt{2}}
\left ( \begin{array}{ccc}
      I & 0 &  -i{\tilde C_0} \\
      0 & \sqrt{2} &  0 \\
      {\tilde C_0} & 0  &  iI
      \end{array} \right ),    \       N=2n+1,
\label{2}
\end{equation}
where $n \times n $ matrix $ {\tilde C_0} $ is like $C_0,$
then  regard the matrix $ D_{\sigma}=DV_{\sigma},$
$V_{\sigma} \in M_{N},$ $(V_{\sigma})_{ik}= \delta_{\sigma_{i},k},$
and $ \sigma \in S(N) $ is a permutation of the $ N$th order.
It is easy to verify that $ D_{\sigma}$
is again a solution of equation (\ref{1}).
Then in symplectic basis the orthogonal Cayley Klein group
$SO(N;j;{\C})$ is described by the matrices
     \begin{equation}
 B_{\sigma}(j)=D_{\sigma}A(j)D^{-1}_{\sigma}
    \label{4}
    \end{equation}
with the additional relations of $ j$-orthogonality
$$
   B_{\sigma}(j)C_{0}B^{t}_{\sigma}(j)=
B^{t}_{\sigma}(j)C_{0}B_{\sigma}(j)=C_{0}.
$$

It should be noted that for  orthogonal groups $(j=1)$
the use of different matrices $D_{\sigma}$ makes no sense because  all
Cartesian coordinates of $O_N$ are equivalent up to a choice of
 its enumerations.
The different situation is for Cayley-Klein groups $(j \not =1).$
Cartesian coordinates
$(1,k)\xi_{k}, \ k=1,\dots,N$ for nilpotent values of some or all
parameters $j_k$ are different elements of the algebra $D_{N-1}(j;\C),$
therefore the same group $SO(N;j;\C)$ may be realized by  matrices
$ B_{\sigma} $  with a {\it different} disposition of nilpotent
generators among their elements.
Namely this fact will provide us with different sets of primitive
elements of Hopf algebra in the case of quantum group.

Matrix elements of $ B_{\sigma}(j) $ are as follows
\begin{equation}
     \begin{array}{l}
(B_{\sigma})_{n+1,n+1}= b_{n+1,n+1},  \\
(B_{\sigma})_{kk}=b_{kk}+i\tilde{b}_{kk}(\sigma_{k},
\sigma_{k'}), \quad
(B_{\sigma})_{k'k'}=b_{kk}-i\tilde{b}_{kk}(\sigma_{k},
\sigma_{k'}), \\
(B_{\sigma})_{kk'}=b_{k'k}-i\tilde{b}_{k'k}(\sigma_{k},
\sigma_{k'}), \quad
(B_{\sigma})_{k'k}=b_{k'k}+i\tilde{b}_{k'k}(\sigma_{k},
\sigma_{k'}), \\
(B_{\sigma})_{k,n+1}=b_{k,n+1}(\sigma_k, \sigma_{n+1}) -
i\tilde{b}_{k,n+1}(\sigma_{n+1}, \sigma_{k'}), \\
(B_{\sigma})_{k',n+1}=b_{k,n+1}(\sigma_k, \sigma_{n+1}) +
i\tilde{b}_{k,n+1}(\sigma_{n+1}, \sigma_{k'}), \\
(B_{\sigma})_{n+1,k}=b_{n+1,k}(\sigma_k, \sigma_{n+1}) +
i\tilde{b}_{n+1,k}(\sigma_{n+1}, \sigma_{k'}), \\
(B_{\sigma})_{n+1,k'}=b_{n+1,k}(\sigma_k, \sigma_{n+1}) -
i\tilde{b}_{n+1,k}(\sigma_{n+1}, \sigma_{k'}), \; k \neq p, \\
(B_{\sigma})_{kp}=b_{kp}(\sigma_k, \sigma_{p})+
b_{kp}'(\sigma_{k'},\sigma_{p'})+
i\tilde{b}_{kp}(\sigma_{k}, \sigma_{p'})-
i\tilde{b}_{kp}'(\sigma_{k'}, \sigma_{p}),  \\
(B_{\sigma})_{kp'}=b_{kp}(\sigma_k, \sigma_{p})-
b_{kp}'(\sigma_{k'},\sigma_{p'})-
i\tilde{b}_{kp}(\sigma_{k}, \sigma_{p'})-
i\tilde{b}_{kp}'(\sigma_{k'}, \sigma_{p}),  \\
(B_{\sigma})_{k'p}=b_{kp}(\sigma_k, \sigma_{p})-
b_{kp}'(\sigma_{k'},\sigma_{p'})+
i\tilde{b}_{kp}(\sigma_{k}, \sigma_{p'})+
i\tilde{b}_{kp}'(\sigma_{k'}, \sigma_{p}),  \\
(B_{\sigma})_{k'p'}=b_{kp}(\sigma_k, \sigma_{p})+
b_{kp}'(\sigma_{k'},\sigma_{p'})-
i\tilde{b}_{kp}(\sigma_{k}, \sigma_{p'})+
i\tilde{b}_{kp}'(\sigma_{k'}, \sigma_{p}) . \\
\end{array}
   \label{5}
\end{equation}
Here $ b,  b', \tilde{b}, \tilde{b}' \in \C$
are expressed by the matrix elements of $A $ with the formula
\begin{displaymath}
     \begin{array}{ll}
b_{n+1,n+1}=a_{\sigma_{n+1},\sigma_{n+1}}, &   \\
b_{n+1,k}=\displaystyle{\frac{1}{\sqrt{2}}a_{\sigma_{n+1},
\sigma_{k}}},
&
b_{k,n+1}=\displaystyle{\frac{1}{\sqrt{2}}a_{\sigma_k,\sigma_{n+1}}},
\\
\tilde{b}_{k,n+1}=
\displaystyle{\frac{1}{\sqrt{2}}a_{\sigma_{k'},\sigma_{n+1}}}, &
\tilde{b}_{n+1,k}=
\displaystyle{\frac{1}{\sqrt{2}}a_{\sigma_{n+1},\sigma_{k'}}}, \\
b_{kk}=\displaystyle{\frac{1}{2}
(a_{\sigma_k\sigma_k}+a_{\sigma_{k'}\sigma_{k'}})}, &
\tilde{b}_{kk}=\displaystyle{\frac{1}{2}
(a_{\sigma_{k}\sigma_{k'}}-a_{\sigma_{k'}\sigma_{k}})}, \\
b_{k'k}=\displaystyle{\frac{1}{2}
(a_{\sigma_k\sigma_k}-a_{\sigma_{k'}\sigma_{k'}})}, &
\tilde{b}_{k'k}=\displaystyle{\frac{1}{2}
(a_{\sigma_k\sigma_{k'}}+a_{\sigma_{k'}\sigma_{k}})}, \\
b_{kp}=\displaystyle{\frac{1}{2}a_{\sigma_k\sigma_p}}, \quad
b'_{kp}=\displaystyle{\frac{1}{2}a_{\sigma_{k'}\sigma_{p'}}}, &
\tilde{b}_{kp}=\displaystyle{\frac{1}{2}a_{\sigma_k\sigma_{p'}}},
\quad
\tilde{b}'_{kp}=\displaystyle{\frac{1}{2}a_{\sigma_{k'}\sigma_p}},
 \;  k \neq p, \\
\end{array}
\end{displaymath}

Let us observe that the elements $b$ of $ B_{\sigma}(j)$
are obtained from the elements $ b^* $ of $ B_{\sigma}(j=1)$
by multiplications on some products of  parameters $j,$ namely
\begin{equation}
     \begin{array}{ll}
b^*_{n+1,n+1}=b_{n+1,n+1}, &
b^*_{kk}=b_{kk}, \quad
b^*_{k'k}=b_{k'k}, \\
\tilde{b}^*_{kk}=(\sigma_k,\sigma_{k'})\tilde{b}_{kk}, &
\tilde{b}^*_{k'k}=(\sigma_k,\sigma_{k'})\tilde{b}_{k'k}, \\
b^*_{k,n+1}=(\sigma_k,\sigma_{n+1})b_{k,n+1}, &
b^*_{n+1,k}=(\sigma_k,\sigma_{n+1})b_{n+1,k},  \\
\tilde{b}^*_{k,n+1}=(\sigma_{k'},\sigma_{n+1})\tilde{b}_{k,n+1}, &
\tilde{b}^*_{n+1,k}=(\sigma_{k'},\sigma_{n+1})\tilde{b}_{n+1,k}, \\
b^*_{kp}=(\sigma_k,\sigma_p)b_{kp}, &
b^{*'}_{kp}=(\sigma_{k'},\sigma_{p'})b'_{kp}, \\
\tilde{b}^*_{kp}=(\sigma_k,\sigma_{p'})\tilde{b}_{kp}, &
\tilde{b}^{*'}_{kp}=(\sigma_{k'},\sigma_p)\tilde{b}'_{kp},
\;  k \neq p.  \\
\end{array}
   \label{6}
\end{equation}
A transformation of group by multiplications of some or all its
group parameters on zero tending parameter $\epsilon$ is  named
as group contraction \cite{IW}, if a new group is obtained in the limit.
The formulas (\ref{6}) are just an example of such transformation,
where the nilpotent values $ j_{k}=\iota_{k} $ are used instead of the
limit $ \epsilon \rightarrow 0. $
In other words group contractions may be described mathematically
correctly by the replacement of real or complex group parameters
with a new one's which are elements of Pimenov algebra
$ {\bf D}_n(\iota; {\C}). $
In our case such replacement is made for matrix elements.

Let us consider as an example the group
 $ SO(3;j;\C).  $ For identical permutation  $ \sigma =(1,2,3) $
the matrix $ D_{\sigma} $ is given by equation (\ref{2}) for $ N=3 $
and in symplectic basis the group $ SO(3;j;\C) $
is described by the matrices
  \begin{displaymath}
        B_{\sigma}(j) =
     \left( \begin{array}{ccc}
 b_{11}+ij_{1}j_{2}\tilde{b}_{11} & j_{1}b_{12}-ij_{2}\tilde{b}_{12}
& b_{31}-ij_{1}j_{2}\tilde{b}_{31} \\
j_{1}b_{21}+ij_{2}\tilde{b}_{21}  & b_{22}
& j_{1}b_{21}-ij_{2}\tilde{b}_{21}       \\
 b_{31}+ij_{1}j_{2}\tilde{b}_{31} & j_{1}b_{12}+ij_{2}\tilde{b}_{12}
& b_{11}-ij_{1}j_{2}\tilde{b}_{11} \\
     \end{array} \right) .
    %\label{36}
 \end{displaymath}
For $ \sigma =(2,1,3)$ one obtain from equation (\ref{4})
 \begin{displaymath}
        B_{\sigma}(j) =
     \left( \begin{array}{ccc}
 b_{11}+ij_{2}\tilde{b}_{11} & j_{1}b_{12}-ij_{1}j_{2}\tilde{b}_{12}
& b_{31}-ij_{2}\tilde{b}_{31} \\
j_{1}b_{21}+ij_{1}j_{2}\tilde{b}_{21}  & b_{22}
& j_{1}b_{21}-ij_{1}j_{2}\tilde{b}_{21}       \\
 b_{31}+ij_{2}\tilde{b}_{31} & j_{1}b_{12}+ij_{1}j_{2}\tilde{b}_{12}
& b_{11}-ij_{2}\tilde{b}_{11} \\
     \end{array} \right) ,
%    \label{37}
 \end{displaymath}
finally the permutation
     $ \sigma =(1,3,2) $
leads to the matrices
 \begin{displaymath}
        B_{\sigma}(j) =
     \left( \begin{array}{ccc}
 b_{11}+ij_{1}\tilde{b}_{11} & j_{1}j_{2}b_{12}-ij_{2}\tilde{b}_{12}
& b_{31}-ij_{1}\tilde{b}_{31} \\
j_{1}j_{2}b_{21}+ij_{2}\tilde{b}_{21}  & b_{22}
& j_{1}j_{2}b_{21}-ij_{2}\tilde{b}_{21}       \\
 b_{31}+ij_{1}\tilde{b}_{31} & j_{1}j_{2}b_{12}+ij_{2}\tilde{b}_{12}
& b_{11}-ij_{1}\tilde{b}_{11} \\
     \end{array} \right) .
%    \label{38}
 \end{displaymath}
The same matrices are corresponded to three
remaining permutations from the group $S(3)$.

For nilpotent values of both parameters
$ j_{1}=\iota_{1}, j_{2}=\iota_{2} $ we have the complex Galilei
group $ G(1+1;\C)=SO(3;\iota;\C), $ which is realized in Cartesian
basis by the matrices
 \begin{displaymath}
        A(\iota) =
     \left( \begin{array}{ccc}
 1  & \iota_{1}a_{12} & \iota_{1}\iota_{2}a_{13} \\
-\iota_{1}a_{12}  & 1 & \iota_{2}a_{23}       \\
 \iota_{1}\iota_{2}a_{31} & -\iota_{2}a_{23} & 1 \\
     \end{array} \right) ,
%    \label{39}
 \end{displaymath}
where $ a_{31}= -a_{13}+a_{12}a_{23}. $
The relations of $ j$-orthogonality (\ref{3})
have been taken into account. Three different realizations of
Galilei group in symplectic description are as follows
 \begin{displaymath}
B_{\sigma}(\iota) =   \left( \begin{array}{ccc}
 1+i\iota_{1}\iota_{2}\tilde{b}_{11}   &
\iota_{1}b_{12}-i\iota_{2}\tilde{b}_{12}
& -i\iota_{1}\iota_{2}\tilde{b}_{31} \\
-\iota_{1}b_{12}-i\iota_{2}\tilde{b}_{12}  & 1
& -\iota_{1}b_{12}+i\iota_{2}\tilde{b}_{12}       \\
 i\iota_{1}\iota_{2}\tilde{b}_{31} & \iota_{1}b_{12}+i\iota_{2}\tilde{b}_{12}
& 1-i\iota_{1}\iota_{2}\tilde{b}_{11} \\
     \end{array} \right) ,
%    \label{40}
 \end{displaymath}
where $ \tilde{b}_{31}=-b_{12}\tilde{b}_{12}, $
   %  $ \sigma =(2,1,3) $
 \begin{displaymath}
        B_{\sigma}(\iota) =
     \left( \begin{array}{ccc}
 1+i\iota_{2}\tilde{b}_{11}   &
\iota_{1}b_{12}-i\iota_{1}\iota_{2}\tilde{b}_{12} & 0 \\
-\iota_{1}b_{12}+i\iota_{1}\iota_{2}\tilde{b}_{21}  & 1
& -\iota_{1}b_{12}-i\iota_{1}\iota_{2}\tilde{b}_{21}       \\
0 & \iota_{1}b_{12}+i\iota_{1}\iota_{2}\tilde{b}_{12}
& 1-i\iota_{2}\tilde{b}_{11} \\
     \end{array} \right) ,
%    \label{41}
 \end{displaymath}
where $ \tilde{b}_{21}=-\tilde{b}_{12} -b_{12}\tilde{b}_{11}, $
    % $ \sigma =(1,3,2) $
 \begin{displaymath}
        B_{\sigma}(\iota) =
     \left( \begin{array}{ccc}
 1+i\iota_{1}\tilde{b}_{11}   &
\iota_{1}\iota_{2}b_{12}-i\iota_{2}\tilde{b}_{12} & 0 \\
\iota_{1}\iota_{2}b_{21}-i\iota_{2}\tilde{b}_{12}  & 1
& \iota_{1}\iota_{2}b_{21}+i\iota_{2}\tilde{b}_{12}       \\
0 & \iota_{1}\iota_{2}b_{12}+i\iota_{2}\tilde{b}_{12}
& 1-i\iota_{1}\tilde{b}_{11} \\
     \end{array} \right) ,
%    \label{42}
 \end{displaymath}
where $ b_{21}=-b_{12} + \tilde{b}_{11}\tilde{b}_{12}. $

\section{ Contractions of quantum orthogonal groups. }
%\label{sec:level2}

%     \subsection{Definition of quantum group $ SO_q(N) $ }
%
%Let us briefly remind the definition of quantum group \cite{Fad-89}.
%The starting point is an algebra $ \C \langle T_{ik} \rangle $
%of noncommutative polynomials of $ N^2 $ variables ---
%a free associate $\C$-algebra with unit generated by
%$ T_{ik},\  i,k=1, \ldots, N. $
%The commutation relations of the generators $ T_{ik} $
%are defined by the equation
%\begin{equation}
%R_qT_1T_2=T_2T_1R_q,
%\label{9}
%\end{equation}
%where  $ T_1=T \otimes I, \ T_2=I \otimes T \in M_{N^2}({\C}
%\langle t_{ik} \rangle), $
%$ T=(T_{ik})^N_{i,k=1} \in M_N({\C} \langle T_{ik} \rangle ) $ is the
%$ N \times N $ matrix with elements from the algebra
%$ \C \langle T_{ik} \rangle $ and $ I $ is the unit matrix in $ M_N({\C}) $.
%The matrix $ T $ satisfies the additional relations of $q$-orthogonality
%\begin{equation}
%TCT^t=T^tCT=C,
%\label{10}
%\end{equation}
%
%{\it Quotient algebra of} $ \C \langle T_{ik} \rangle $
%{\it by relations} (\ref{9}), (\ref{10}) {\it is called the quantum
%orthogonal group (or functions algebra on quantum group)}
%$ SO_{q}(N) $ and is a Hopf algebra with the following coproduct
% $ \Delta, $ counits $ \epsilon $ and antipode $ S $
%$$
%\Delta T=T \dot {\otimes} T, \quad
%\Delta T_{ik}=\sum^N_{m=1} T_{im} \otimes T_{mk},
%$$
%$$
%\epsilon(T)=I, \; \epsilon(T_{ik})=\delta_{ik},
%$$
%$$
%S(T)=CT^tC^{-1}, \quad  S(T_{ik})=q^{\rho_{i'}-\rho_{k'}}T_{k'i'}, \;
%i,k=1, \ldots, N.
%$$

\subsection{Formal definition of  quantum group $ SO_v(N;j;\sigma) $ }

In the definition of the quantum group $ SO_v(N;j;\sigma)$
we shall follow  \cite{Fad-89}, but
 start with an algebra $ {\bf D} \langle (T_{\sigma})_{ik} \rangle $
of noncommutative polynomials of  $ N^2 $ variables, which are an elements of
the direct product $ {\bf D}_{N-1}(j)\otimes {\C}\langle t_{ik} \rangle. $
More precisely, the elements $  (T_{\sigma})_{ik} $ are obtained
from the elements $ \left( B_{\sigma}(j) \right)_{ik}$ of equations~(\ref{5})
by the replacement of commutative variables $ b,b',\tilde{b},\tilde{b}' $
with the noncommutative variables $ t,t',\tau,\tau',$ respectively.
It is clear that generators $ t,t',\tau,\tau'$ are connected with
the corresponding generators $ t^*,t^{*}{'},\tau^*,\tau^{*}{'}$ of
$SO_q(N)$ in just the same way~(\ref{6}) as elements
$ b,b',\tilde b,\tilde b'$ are connected with
$ b^*,b^{*}{'},\tilde b^*, \tilde b^{*}{'}$.
%The last ones are transformed by equations~(\ref{6}).
One introduces additionally the transformation of the deformation parameters
$q=e^z $ as follows:
\begin{equation}
z=Jv,
\label{11}
\end{equation}
where $ v $ is a new deformation parameter and  $ J $ is some
product of  parameters $ j $ for the present unknown.
Nondegenerate low triangular matrix $ R_{q} \in M_{N^2}(\C) $
is given by
$$
R_q=q\sum^N_{k=1,k \neq k'}e_{kk} \otimes e_{kk} +
$$
$$
+\sum^N_{k,r=1, k \neq r,r'}e_{kk} \otimes e_{rr}+
q^{-1}\sum^N_{k=1,k \neq k'}e_{k'k'} \otimes e_{kk}+
(q-q^{-1})\sum^N_{k,r=1, \ k>r}e_{kr} \otimes e_{rk}-
$$
$$
%\begin{equation}
-(q-q^{-1})\sum^N_{k,r=1, \ k>r}q^{\rho_k-\rho_r}e_{kr} \otimes e_{k'r'}+
e_{pp} \otimes e_{pp},
%\label{7}
%\end{equation}
$$
where the last term is present only for $ N=2n+1 $ and
$ p=(N+1)/2 $. Here $ e_{ik} \in M_N({\C}) $ are the matrix units
$ (e_{ik})_{sm}=\delta_{is}\delta_{km}, $ $ k'=N+1-k, \ r'=N+1-r $ and
\begin{equation}
(\rho_1, \ldots, \rho_N)=
\left \{ \begin{array}{ccc}
     (n-\frac{1}{2}, n-\frac{3}{2}, \ldots , \frac{1}{2},0,-\frac{1}{2},
     \ldots , -n+\frac{1}{2}), \; N=2n+1 \\
     (n-1, n-2, \ldots, 1,0,0,-1, \ldots, -n+1), \; N=2n.
     \end{array} \right.
\label{8}
\end{equation}
Matrix $C$ is as follows
 $$
C=C_0q^{\rho},\; \rho=diag(\rho_1, \ldots, \rho_N), \;
(C_0)_{ik}=\delta_{i'k}, \  i,k=1, \ldots, N,
$$
$$ (C)_{ik}=q^{\rho_{i'}}\delta_{i'k}, \;\;
 (C^{-1})_{ik}=q^{-\rho_i}\delta_{i'k}.
$$
Let $ \tilde R_v(j), C(j) $ be matrices which are obtained from
$R_q, C$
by the replecement of deformation parameter $ z $ with $ Jv: $
$$
R_{v}(j)=R_{q}(z \rightarrow Jv), \quad
C(j)=C(z \rightarrow Jv).
$$
The commutation relations of the generators
$ T_{\sigma}(j) $ are defined by
\begin{equation}
 R_v(j)T_1(j)T_2(j)=T_2(j)T_1(j)R_v(j),
\label{12}
\end{equation}
where $ T_1(j)=T_{\sigma}(j) \otimes I, \
T_2(j)=I \otimes T_{\sigma}(j)   $
and  the additional relations of $ (v,j)$-orthogonality
     \begin{equation}
   T_{\sigma}(j)C(j)T^{t}_{\sigma}(j) =
T^{t}_{\sigma}(j)C(j)T_{\sigma}(j) = C(j).
    \label{13}
    \end{equation}
are imposed.

{\it One defines the quantum orthogonal Cayley-Klein group}
$ SO_{v}(N;j;\sigma) $ {\it as the quotient algebra of}
 $ {\bf D} \langle (T_{\sigma})_{ik} \rangle $
{\it by relations} (\ref{12}),(\ref{13}).
Formally $ SO_{v}(N;j;\sigma)$ is a Hopf algebra with the following coproduct
$ \Delta, $ counit $ \epsilon $ and antipode $ S: $
$$
\Delta T_{\sigma}(j)=T_{\sigma}(j) \dot {\otimes}T_{\sigma}(j),
\quad \epsilon (T_{\sigma}(j))=I, \quad
S(T_{\sigma}(j))=C(j)T^{t}_{\sigma}(j)C^{-1}(j).
$$

In terms of generators $ t,\tau $ the explicit form of antipode is given
in Appendix A, of coproduct is given in Appendix B and
of  $ (v,j)$-orthogonality relations are given in Appendix C.
As far as only secondary diagonal elements of the matrix
 $ C(j) $   are different from zero and for
$ q=1, j=1 $ it is equal to $ C_{0}, $ then
we have the symplectic description of $ SO_v(N;j;\sigma). $

\subsection{Allowed contractions of $ SO_v(N;j;\sigma) $ }

The formal definition of the quantum group
$ SO_{v}(N;j;\sigma) $ should be a real definition of quantum group,
if the proposed construction is a consistent Hopf algebra structure
under nilpotent values of some or all parameters $j$.
Counit $ \epsilon(t_{n+1,n+1})=1, \; \epsilon(t_{kk})=1, \;
k=1,\ldots,n,$ and $\epsilon(t)=\epsilon(\tau)=0 $
for the rest generators
do not restrict the values of  parameters $ j. $
Parameters $j$ are arranged in the expressions for coproduct
$ \Delta $ (Appendix B) exactly as in matrix product of
$ B_{\sigma}(j), $ and as far as the last ones form the group
$ SO(N;j;\C)$ for any values of $j$,  no restrictions
follow from the coproduct. Different situation is with the
antipode $ S $ (Appendix A). Really, for elements
\begin{equation}
(T_{\sigma})_{k'k}=t_{k'k} + i\tau_{k'k}(\sigma_k,\sigma_{k'}),
  \quad  k=1,\ldots,n,
\label{14}
\end{equation}
 antipode is obtained as
\begin{equation}
S((T_{\sigma})_{k'k})=(T_{\sigma})_{k'k}\cdot e^{2J\rho_kv}, \quad
\label{15}
\end{equation}
and depend both on $ \rho_k$ and for the present undetermined
factor $J.$ Antipode is antihomomorphism of Hopf algebra and
therefore has to transform the matrix $ T_{\sigma}(j) $ to a matrix
with the same distribution of the nilpotent parameters $j$ in its elements,
i.e. the right and the left parts of equation (\ref{15})
must be identical elements of
$ {\bf D}_{N-1}(j)\otimes {\C}\langle t_{ik} \rangle. $
For $J=1$ this condition  holds for any values of the parameters $j.$
The case $ J \neq 1$  requires additional discussion.

Next condition which must be taken into account is the
$(v,j)$-orthogonality relations (\ref{13}) (see Appendix C).
In general, for nilpotent values of parameters   $ j_k $
the number of equations (\ref{13}) are increased as compared
with the case $ j_k=1,\ k=1,\ldots,n$ because it is necessary equate to
each other terms with nilpotent generators and their products independently.
 Then the number of contracted quantum
group generators are decreased as
compared with the initial $SO_q(N)$.
For example, for $ j_k=\iota_k, \; k=1,\ldots,n, \; J=1 $
the only nonzero generators of quantum group $ SO_q(2n+1;\iota;\sigma_{k}=k) $
are $ \tau_{kk}, \tau_{k'k}, t_{kk}=1, \; k=1,\ldots,n, t_{n+1,n+1}=1, $
i.e. as a result of contraction we have the Hopf algebra with the number
of generators equal to $ 2n-1$ which is less then $ N(N-1)/2=n(2n+1). $

On the other hand most interesting are such contractions,
when the number of generators is conserved. It is necessery
for this that the number of equations in
$ (v,j)$-orthogonality relations is not changed as compared
with the initial quantum group. It is possible when nilpotent generators
 appeare in equation (\ref{13}) either with the powers greater or
equal two (and then the corresponding terms are equal to zero) or
as homogeneous multipliers.
Taking into account all these arguments and using the explicit expressions
for antipode and $(v,j)$-orthogonality we can find
possible contractions of quantum orthogonal groups, which are
described by the following theorems.

{ \bf Theorem 1.}   { \it
If the deformation parameter is not transformed $J=1$,
then the following maximal $n$-dimensional contraction of the
orthogonal quantum group  $ SO_v(N;j;\sigma), \ N=2n+1 $ is allowed:
 \begin{equation}
 j_{2s}=\iota_{2s}, \; s=1,\ldots,m, \;
j_{2r+1}=\iota_{2r+1},  \; r=m,\ldots,n-1, \; 0 \leq m \leq n,
  \label{16}
  \end{equation}
for example, for permutation $\sigma $:
$
\sigma_{n+1}=2m+1, \;
\sigma_{s}=2s-1, \; \sigma_{s'}=2s, \; s=1,\ldots,m, \;
\sigma_{r}=2r, \; \sigma_{r'}=2r+1, \; r=m+1,\ldots,n.
$ }

{ \it Proof.} From the explicit form of $ (v,j)$-orthogonality
(Appendix C) it follows that if all multipliers
$ (\sigma_{k},\sigma_{k'}), \; k=1,\ldots,n $
are equal to one (or equivalently $ \bigcup_{k=1}^{n}(\sigma_{k},\sigma_{k'})=1$),
then under conditions of theorem all products of the
parameters $ j $ in (\ref{13}) are equal to zero, otherwise are appeared
in these equations as homogeneous multipliers.
\rule{1ex}{1ex}

{ \bf Theorem 2.}   { \it
If the deformation parameter is not transformed $J=1$,
then the following maximal n-dimensional contraction of the quantum orthogonal
group  $ SO_v(N;j;\sigma), \ N=2n $ is allowed:
\begin{eqnarray}
& & j_{2s}=\iota_{2s}, \; s=1,\ldots,m-1, \;
j_{2p-1}=\iota_{2p-1},  \; p=m,\ldots,u, \; \nonumber \\
& & j_{2r}=\iota_{2r},  \; r=u,\ldots,n-1, \;
1 \leq m \leq u \leq n,
  \label{17}
  \end{eqnarray}
for example, for permutation $ \sigma $:
$ \sigma_{n}=2m-1, \;\sigma_{n'}=2u, \;
\sigma_{s}=2s-1, \; \sigma_{s'}=2s, \; s=1,\ldots,m-1, \;
\sigma_{p}=2p, \; \sigma_{p'}=2p+1, \; p=m,\ldots,u-1, \;
\sigma_{r}=2r+1, \; \sigma_{r'}=2r, \; r=u,\ldots,n-1.
$ }

{ \it Proof.} Similar to the proof theorem 1, except $k=1,\ldots,n-1$.
\rule{1ex}{1ex}

{\bf Remark 1.}
It should be noted that as $ \sigma $ may be taken any permutation
with the properties $ (\sigma_{k},\sigma_{k'})=1, \; k=1,\ldots,n $
(or $ n-1$).

{\bf Remark 2.}
Admissible contractions for number of parameters  $ j_k $
less then $n$ are obtained from (\ref{16}),(\ref{17}),
by setting part of $ j_{2s},  j_{2p-1},  j_{2r},$ $ j_{2r+1} $ equal to one.

We return to the antipode  (\ref{15}) for $ J \neq 1.$
As far as $ \rho_{n+1}=0 $ for $ N=2n+1,$ and $ \rho_{n} = \rho_{n'}=0 $
for $N=2n,$ (\ref{8}) we shall regard these two cases separately.

{ \bf Theorem 3.}   { \it
If the deformation parameter is transformed $(J \not =1)$,
then the following contractions of the quantum orthogonal
group  $ SO_v(N;j;\sigma),$ $ N=2n+1 $ are allowed:

1. For  $ J=j_{n+1}, $

\noindent
\begin{tabular}{ll}
a) $ j_{n+1}=\iota_{n+1},$ & if
$ 1<\sigma_{n+1}<n+1;  $ \\
b) $ j_{n+1}=\iota_{n+1}, \; j_{1}=1,\iota_{1},$
& if $ \sigma_{n+1} = 1.$ \\
\end{tabular}

2. For  $ J=j_{n}, $

\noindent
\begin{tabular}{ll}
a) $ j_{n}=\iota_{n}, $ & if
 $n+1<\sigma_{n+1}<2n+1;$ \\
b) $ j_{n}=\iota_{n}, \;  j_{2n}=1,\iota_{2n},$
& if $ \sigma_{n+1} =2n+1. $ \\
\end{tabular}

3. For  $  J=j_nj_{n+1},$

\noindent
\begin{tabular}{ll}
$j_{n}=1,\iota_{n},\; j_{n+1}=1,\iota_{n+1}, $
& if  $ \sigma_{n+1}=n+1. $
\end{tabular}
}

{ \it Proof.  }
If $ J \sim \iota, $ then $ e^{J\rho v}=1+J\rho v, $
and equation (\ref{15}) is rewritten as
\begin{equation}
S(T_{k'k})=t_{k'k} + i\tau_{k'k}(\sigma_k,\sigma_{k'}) +
2t_{k'k} \rho_kvJ+2i\tau_{k',k}(\sigma_{k},\sigma_{k'})\rho_kvJ.
\label{18}
\end{equation}
The terms with factor $J$ may be added only with the terms
with factors $(\sigma_k,\sigma_{k'}), \ k=1,\dots,n, $ therefore
\begin{equation}
J=\bigcap_{k=1}^{n}(\sigma_{k},\sigma_{k'}),
\label{19}
\end{equation}
i.e. multiplier $J$ is the product of all nilpotent generators
of Pimenov algebra, which are simultaneously contained in all
$(\sigma_k,\sigma_{k'}),\ k=1,\dots,n.$
Intersection is not empty if all $ \sigma_{k}  $ are less then all
$\sigma_{p'},$ i.e.
$\sigma_{k} < \sigma_{p'}, \; \forall k,p=1, \ldots ,n. $
The terms $ J(\sigma_{k},\sigma_{k'}) $ are equal to zero.
Taking into account that $ \sigma_{n+1} $ is not contained
in equation (\ref{19}) we may find unknown multiplier $ J. $

A. If $ \sigma_{n+1}<n+1, $ then
$ \max \sigma_p=n+1 $ and
$ \min  \sigma_{p'}=n+2, $ therefore
$ J=\bigcap_{p=1}^{n}(\sigma_{k},\sigma_{k'})=(n+1,n+2)=j_{n+1}. $

B. If  $ \sigma_{n+1}>n+1, $ then
$ \max  \sigma_p=n $ and
$ \min  \sigma_{p'}=n+1, $ therefore
$ J=\bigcap_{p=1}^{n}(\sigma_{k},\sigma_{k'})=(n,n+1)=j_{n}. $

C. If  $ \sigma_{n+1}=n+1, $ then
$ \max  \sigma_p=n $ and
$ \min  \sigma_{p'}=n+2, $ therefore
$ J=\bigcap_{p=1}^{n}(\sigma_{k},\sigma_{k'})=(n,n+2)=j_nj_{n+1}. $

Let us return to (\ref{18}) with regard of obtained possible values $ J. $

1. For $ J=j_{n+1}. $

a) For permutations $\sigma$ with
$ 1 < \sigma_{n+1} < n+1 $
only one contraction $ j_{n+1}=\iota_{n+1}$ is allowed.

b) If $ \sigma_{n+1} = 1, $ then all products
$ (\sigma_k, \sigma_{k'}), \; k=1,\ldots,n $
do not contain parameter $ j_{1}, $ therefore it is not appeared in
$ (T_{\sigma})_{kk'}, \; S((T_{\sigma})_{kk'}). $
Consequently for nilpotent value of $ j_{1} $
the above mentioned matrix elements and their antipodes are the same
elements of $ {\bf D}_{N-1}(j)\otimes \C\langle t_{ik} \rangle,$
i.e. for permutations $\sigma$ with $ \sigma_{n+1} = 1 $
two dimensional contraction
$ j_{n+1}=\iota_{n+1}, \; j_{1}=\iota_{1} $ is allowed.

2. For $ J=j_{n}. $

a) For permutations $\sigma$ with
$ n+1 < \sigma_{n+1} < 2n+1 $
only one contraction $ j_{n}=\iota_{n}$ is allowed.

b) For permutations $\sigma$ with
$ \sigma_{n+1} =2n+1 $
two dimensional contraction
$ j_{n}=\iota_{n}, \;  j_{2n}= \iota_{2n}$ is allowed since
all products $ (\sigma_k, \sigma_{k'}), \; k=1,\ldots,n $
do not contain parameters $ j_{2n}. $

3. For  $ J=j_{n}j_{n+1} $
the permutations $\sigma$ with $  \sigma_{n+1}=n+1$
are regarded and both parameters $ j_n,j_{n+1}$ may be independently
equal to nilpotent values, therefore one have three contractions:
 $ j_{n}=\iota_{n},\; j_{n+1}=1; \quad
   j_{n}=1,\; j_{n+1}=\iota_{n+1} \ and \
   j_{n}=\iota_{n},\; j_{n+1}=\iota_{n+1}.  $

We have found the admissible contractions by analysis of
antipode of the matrix elements
$ (T_{\sigma})_{k'k},\;  k=1,\ldots,n. $
One may verify that the antipode of remaining elements of
$ T_{\sigma} $ leads to the same admissible contractions.
In this sense the selected elements $ \left( T_{\sigma}\right )_{k'k} $
are most informative.

Using the explicit form of $ (v,j)$-orthogonality (Appendix C),
it is easy to verify that under the conditions of theorem all
products of the parameters  $j$ are equal to one or zero, otherwise are
appeared in (\ref{13}) as homogeneous multipliers.
\rule{1ex}{1ex}

{ \bf Theorem 4.}   { \it
     If the deformation parameter is transformed $(J\neq1)$,
then the following contractions of the quantum orthogonal group
$SO_v(N;j;\sigma),\ N=2n $ are allowed:

1. For $ J=j_{n}, $

\noindent
\begin{tabular}{ll}
a) $ j_{n}=\iota_n, $ & if
$ \sigma_{n}>1,\;  \sigma_{n'}<2n; $ \\
b) $ j_{n}=\iota_n, \; j_{1}=1,\iota_{1}, $ & if
$ \sigma_{n}=1,\;  \sigma_{n'}< 2n; $ \\
c) $ j_{n}=\iota_n, \; j_{2n-1}=1,\iota_{2n-1}, $ & if
$\sigma_n > 1,\; \sigma_{n'}= 2n; $\\
d) $ j_{n}=\iota_n, \; j_{1}=1,\iota_{1},
\; j_{2n-1}=1,\iota_{2n-1}, $ &
if $\sigma_n=1,\;  \sigma_{n'}= 2n.  $
\end{tabular}

2. For $ J=j_{n-1}. $

\noindent
\begin{tabular}{ll}
a) $ j_{n-1}=\iota_{n-1}, $ & if
$ \sigma_{n'}< 2n; $\\
b) $ j_{n-1}=\iota_{n-1}, \; j_{2n-1}=1,\iota_{2n-1}, $
& if $ \sigma_n<2n-1, \; \sigma_{n'}=2n; $\\
c) $ j_{n-1}=\iota_{n-1}, \; j_{2n-2}=1,\iota_{2n-2}, $&\\
$ \quad \ j_{2n-1}=1,\iota_{2n-1},  $
& if $\sigma_n=2n-1,\;  \sigma_{n'}= 2n.  $
\end{tabular}

3. For  $ J=j_{n+1}. $

\noindent
\begin{tabular}{ll}
a) $ j_{n+1}=\iota_{n+1}, $
& if $ \sigma_{n}>1; $ \\
b) $ j_{n+1}=\iota_{n+1}, \; j_{1}=1,\iota_{1}, $
& if $ \sigma_{n}=1, \; \sigma_{n'}>2; $ \\
c) $ j_{n+1}=\iota_{n+1}, \; j_{1}=1,\iota_{1}, \;
j_{2}=1,\iota_{2},  $ & if $ \sigma_{n}=1,\;  \sigma_{n'}= 2.  $
\end{tabular}

4. For  $ J=j_{n-1}j_{n}. $

\noindent
\begin{tabular}{ll}
a) $ j_{n-1}=\iota_{n-1}, \; j_{n}=\iota_{n}, $
& if $ \sigma_{n'}<2n; $ \\
b) $ j_{n-1}=\iota_{n-1}, \; j_{n}=\iota_n, \; j_{2n-1}=1,\iota_{2n-1}, $
& if $ \sigma_{n'}=2n.  $
\end{tabular}

5. For  $ J=j_{n}j_{n+1}. $

\noindent
\begin{tabular}{ll}
a) $ j_{n}=\iota_{n}, \; j_{n+1}=\iota_{n+1}, $
& if $\sigma_{n}>1; $ \\
b) $ j_{1}=1,\iota_{1}, \; j_{n}=\iota_{n}, \; j_{n+1}=\iota_{n+1}, $
& if $ \sigma_{n}=1. $
\end{tabular}

6. For  $ J=j_{n-1}j_{n}j_{n+1}. $

\noindent
\begin{tabular}{ll}
a) $ j_{n-1}=1,\iota_{n-1}, \; j_{n}=1,\iota_{n}, \; j_{n+1}=1,\iota_{n+1}, $
& if
 $ \sigma_{n}=n,\; \sigma_{n'}=n+1. $
\end{tabular}
    }

 { \it Proof.} The matrix elements $ (T_{\sigma})_{k'k} $
and their antipodes are described by equations
(\ref{14}),(\ref{15}) with $ k=1, \ldots,n-1. $
(As far as $ \rho_{n}=0, $ then $ S((T_{\sigma})_{nn'})= S((T_{\sigma})_{n,n+1})
 =(T_{\sigma})_{n,n+1})$.
Therefore $ J $ is given by equations~(\ref{19}) with the replacement of $ n $
by $ n-1. $ Considering $ \sigma_{k} < \sigma_{p'}, \; \forall k,p=1, \ldots n,$
one find admissible values of $J.$

A. Let $ \sigma_n<n $ and $\sigma_{n'}>n+1, $ then
$ \max  \sigma_k=n $ and
$ \min  \sigma_{k'}=n+1, $ hence $ J=(n,n+1)=j_n. $

B. Let  $ \sigma_{n'}>\sigma_n>n, $ then
$ \max  \sigma_p=n-1 $ and
$ \min  \sigma_{p'}=n,$
hence $ J=(n-1,n)=j_{n-1}. $

C. Let $ n\geq \sigma_{n'}>\sigma_n, $ then
$ \max  \sigma_p=n+1 $ and
$ \min  \sigma_{p'}=n+2, $
hence  $  J=(n+1,n+2)=j_{n+1}. $

D. Let $ \sigma_{n}=n, \;   \sigma_{n'}>n+1, $ then
$ \max  \sigma_p=n-1 $ and
$ \min  \sigma_{p'}=n+1, $
hence   $ J=(n-1,n+1)=j_{n-1}j_n. $

E. Let $ \sigma_{n}<n, \;   \sigma_{n'}=n+1, $ then
$ \max  \sigma_p=n $ and
$ \min  \sigma_{p'}=n+2, $
hence  $ J=(n,n+2)=j_{n}j_{n+1}. $

F. Let $ \sigma_{n}=n, \; \ \sigma_{n'}=n+1, $ then
$ \max  \sigma_p=n-1 $ and
$ \min  \sigma_{p'}=n+2,$
hence   $ J=(n-1,n+2)=j_{n-1}j_{n}j_{n+1}. $

The analysis of equations~(\ref{15}), with due regard for obtained
possible values of $J,$ leads to the admissible contractions
of the theorem. Using the explicit form of $ (v,j)$-orthogonality (Appendix C),
it is easy to verify that under the conditions of theorem all products of
the parameters  $ j $ are equal to one or zero, otherwise are appeared in
(\ref{13}) as homogeneous multipliers.
\rule{1ex}{1ex}

Hopf algebra  $ SO_q(N;j;\sigma),\ N=2n+1$ has $n$
primitive elements which correspond to $n$ nonintersecting
$2\times2$ submatricies of the Cartesian matrix $ A(j)$
composed from elements $ a_{\sigma_{k}\sigma_{k}}, a_{\sigma_{k}\sigma_{k'}},
a_{\sigma_{k'}\sigma_{k}}, a_{\sigma_{k'}\sigma_{k'}}, k=1,\ldots,n. $
Under the transition to the symplectic basis they are transformed to $n$
diagonal $ 2\times2 $ submatricies
diag$((B_{\sigma})_{kk},(B_{\sigma})_{k'k'}) =
$ diag$(b_{kk}+i\tilde{b}_{kk}(\sigma_{k},\sigma_{k'}),
     b_{kk}-i\tilde{b}_{kk}(\sigma_{k},\sigma_{k'})), \,
     k=1,\ldots,n, $
see (\ref{5}). Each such matrix is either one parameter
rotation subgroup $ SO(2),$ if $  (\sigma_{k},\sigma_{k'})=1, $
or one parameter Galilei transformation $ SO(2;j=\iota)=G(1,1)$,
if  $ (\sigma_{k},\sigma_{k'})=\iota.$
Therefore, if the deformation parameter $z$ is fixed
$(J=1)$, then {\it all} primitive elements
of the contracted quantum orthogonal groups  correspond to
Euclidean rotation  $ SO(2).$
If the deformation parameter is transformed $ z=\iota v$,
then { \it all } primitive elements  correspond
to Galilei transformation  $ SO(2;j=\iota)=G(1,1).$ The same is true for the
contracted quantum groups $ SO_q(N;j;\sigma),\ N=2n.$
Let us note that contractions of quantum orthogonal algebras
with  different sets of primitive elements have been discussed in
\cite{VG},\cite{T}.

Quantum orthogonal groups have contractions with the
same nilpotent parameters $j$ both with a fixed deformation
parameter and with a transformed one. For example, the quantum group
$ SO_q(2n+1;j;\sigma)$ for even $ n=2p $  at
$ \sigma_{n+1}=1 $ according to (\ref{16}) has contraction
$ j_{n}=\iota_{n}, \; j_{n+1}=\iota_{n+1}, \; J=1$
and according to  3 of Theorem 3 has the same two-dimensional contraction,
but $  J=\iota_{n}\iota_{n+1}. $ Quantum group
$ SO_q(2n;j;\sigma)$ for odd $ n=2p-1 $ at
$ \sigma_{n}=n, \; \sigma_{n'}=n+1 $ according to (\ref{17}) has
contraction $ j_{n-1}=\iota_{n-1}, j_{n}=\iota_{n}, \; j_{n+1}=\iota_{n+1},
\; J=1$ and according to 6 of Theorem 4 has the same three-dimensional
contraction but $ J=\iota_{n-1}\iota_{n}\iota_{n+1}.$
Let us stress that the cases $ J=1$ and $J \sim  \iota $
are realized for {\it different} sets of primitive elements in
Hopf algebras $ SO_q(4p+1; \iota_{n}, \iota_{n+1}; \sigma)$ and
$ SO_q(4p-2; \iota_{n-1}, \iota_{n}, \iota_{n+1}; \sigma), $ respectively.

Let permutation $\sigma$ be identical, i.e.
$\sigma_k=k,$ $\sigma_{k'}=k',$ $\sigma_{n+1}=n+1.$
It follows from theorems  1 and 2 that there are no contractions of
quantum orthogonal group $SO_q(N;j)$ with fixed deformation parameter $(J=1)$.
For $N=2n+1$ from theorem 3 we obtain three possible contractions
$j_n=1,\iota_n,$ $j_{n+1}=1,\iota_{n+1}$
(both parameters $j_n$ and $j_{n+1}$ independently take nilpotent values)
and deformation parameters is transformed by (\ref{11}), with $J=j_nj_{n+1}.$
For $N=2n$ from theorem 4 we obtain seven admissible contractions:
$j_{n-1}=1, \iota_{n-1},$  $j_{n}=1, \iota_{n},$
$j_{n+1}=1, \iota_{n+1},$ where deformation parameter is multiplied by
$J=j_{n-1}j_nj_{n+1}.$ It should be considered in papers
\cite{GKK},\cite{Sb-97}
just these allowed contractions.

\section{Quantum complex kinematic groups}
%\label{sec:level3}

Kinematic groups are  motion groups of the maximal homogeneous
four-dimensional (one time and three space coordinates)
space--time models \cite{BLL}. All these groups may be obtained from the
real orthogonal group $ SO(5;\R) $ by contractions and analytic
continuations \cite{2}.
If one introduce Beltrami coordinates
     $ \xi_{k}=x_{k+1}/x_{1}, \; k=1,2,3,4 $
and one interpret $\xi_1$ as a time axis while the rest three --
as a space axes, then Galilei group $ G(1,3)=SO(5;\iota_{1},\iota_{2},1,1) $
is the motion group of the nonrelativistic space-time with zero curvature,
Newton groups $ N^{\pm}(1,3)=SO(5;j_{1}=1,i;\iota_{2},1,1) $
are the motion groups of the nonrelativistic space-time with positive
and negative curvature, respectively. Poincare group
$ P(1,3)=SO(5;\iota_{1},i,1,1) $ is the motion group of the relativistic
space-time with zero curvature and $ S^{\pm}(1,3)=SO(5;j_{1}=1,i;i,1,1) $
are the motion groups of the anti de Sitter space-time (positive curvature)
and de Sitter space-time (negative curvature).

If one interpret three first Beltrami coordinates as a space axes while the
last one as a time axis, then the three exotic Carroll kinematics are
obtained, namely $ C^{0}(1,3)=SO(5;\iota_{1},1,1,\iota_{4}),$
with zero curvature, $ C^{\pm}(1,3)=SO(5;j_{1}=1,i;1,1,\iota_{4}),$
with positive and negative curvature.
%The Carroll space and time are as it were
%interchanged their properties as compared with Galilei kinematic.
%The Galilei time is absolute, i.e. two simultaneous in some reference
%frame events remain simultaneous in any reference frame which is obtained
%by Galilei boost (or space-time rotation) from the initial one.
%On the contrary in Carroll kinematic the space is absolute, i.e. two
%events with the equal spatial coordinates in any reference frame have
%the same spatial coordinates in any reference frame which is obtained
%from the initial one by space-time rotation.

The groups $ N^{\pm}(1,3)$ are the real forms of the complex Newton
group $ N(4),$ Poincare group $P(1,3)$ is the real form of the
complex Euclid group $E(4),$ the groups $ C^{\pm}(1,3) $
are the real forms of the complex Carroll group $C(4).$
In this paper the quantum deformations of the complex orthogonal groups
are regarded, therefore whith the help of contractions of
$SO_q(5)$ the quantum analogs
of the complex kinematic groups may be obtained. Possible contractions
of the complex quantum group $SO_{q}(5;j;\sigma)$ are
described by the theorems 1,3 and are
as follows:
for $ J=1,$
$ j_{1}=1,\iota_{1},\; j_{3}=1,\iota_{3} $ with
$ \sigma=(2,4,1,5,3); \; $
$ j_{2}=1,\iota_{2},\; j_{3}=1,\iota_{3} $ with
$ \sigma=(1,4,3,5,2); \; $
$ j_{2}=1,\iota_{2},\; j_{4}=1,\iota_{4} $ with
$ \sigma=(1,3,5,4,2); \; $
for $ J=\iota_{2}, $
$ j_{2}=\iota_{2},\; j_{4}=1,\iota_{4} $ with
$ \sigma=(1,2,5,3,4); \; $
for $ J=\iota_{3}, $
$ j_{1}=1,\iota_{1},\; j_{3}=\iota_{3} $ with
$ \sigma=(2,3,1,4,5); \; $
for $ J=\iota_{2}\iota_{3}, $ $ j_{2}=\iota_{2},\; j_{3}=\iota_{3} $ with
$ \sigma=(1,2,3,4,5). $
Thus if  deformation parameter remains unchanged $( J=1),$ then
we have the quantum analog of  Euclidean group $ E_{q}(4) $
for  $ j_{1}=\iota_{1}, j_{2}=j_{3}=j_{4}=1,\; \sigma=(2,4,1,5,3); $
of  Newton group $ N_{q}(4) $
for $ j_{2}=\iota_{2}, j_{1}=j_{3}=j_{4}=1,\; \sigma=(1,4,3,5,2)$
and of  Carroll group $ C_{q}(4) $
for $ j_{4}=\iota_{4}, j_{1}=j_{2}=j_{3}=1,\; \sigma=(1,3,5,4,2).$
If  deformation parameter is transformed under contraction $ z=\iota_{2}v,$
then we have one more quantum deformation of  Newton group $ N_{v}(4) $
for $ j_{2}=\iota_{2}, j_{1}=j_{3}=j_{4}=1,\; \sigma=(1,2,5,3,4), $
which is not isomorphic to the previous one. Two primitive elements of
$ N_{q}(4)$ correspond to the elliptic translation along the time axis
$ t $ and to the rotation in the space plane $ \{r_{2},r_{3}\} $ (both
are isomorphic to $ SO(2)$), while primitive elements of $ N_{v}(4)$
correspond to the flat translation along the spatial axis $ r_{2}$
and to  Galilei boost in the space-time plane $ \{t,r_{1}\}$
(both are isomorphic to Galilei group $ SO(2;j_{2}=\iota_{2})=G(1,1)).$
We did not obtain the quantum deformations of the complex Galilei
$G(4)$ and Carroll $C^0(4)$ groups.

According to correspondence principle a new physical theory must include
an old one as a particular case. For space-time theory this principle is
realized as the chain of limit transitions: general relativity passes to
special relativity, when space-time curvature tends to zero, and special
relativity passes to classical physics, when light velosity tends to infinity.
For kinematical groups this corresponds to the chain of contractions:
  \begin{equation}
  S^{\pm}(1,3)\stackrel{ K \rightarrow 0}{\longrightarrow}
  P(1,3)\stackrel{c \rightarrow \infty }{\longrightarrow}G(1,3).
    \label{20}
     \end{equation}
As it was mentioned above there is no quantum deformation of the complex
Galilei group in our scheme, therefore we are not able to construct
the standard quantum
analog of the full chain of contractions (\ref{20}), even at the level
of complex groups.
This means that (at least standard) quantum
deformation of the flat nonrelativistic (1+3) space-time  does not exist
in Cayley--Klein scheme.

\section{Conclusion}

From the contraction viewpoint  Hopf algebra structure of quantum
orthogonal group is more rigid as compared with a group one.
Cayley-Klein groups are obtained \cite{2} from $SO(N;j)$
for all nilpotent values
of  parameters $j_k,\ k=1,\ldots,N-1,$  whereas their quantum
deformations  exist only for  some of them $(\leq[{N \over 2}])$.
The main restrictions on contractions are appeared from antipode
(\ref{15}). In particular, contractions of quantum orthogonal groups with 
transformed deformation parameter $z=Jv, \ J\neq 1$ are possible
only due to some parameters (\ref{8}), which characterize the matrix $R_q$,
are equal to zero, namely $\rho_{n+1}=0$ for $N=2n+1$ and
$\rho_n=\rho_{n'}=0$ for $N=2n$. In this sence such contractions
are exclusive and complementary to contractions
with untransformed deformation parameter.

It should be noted that among the contracted for equal number of
parameters $j$ quantum orthogonal groups may be isomorphic, as Hopf
algebra quantum groups. Quantum groups isomorphism is not regarded
in this paper.

Unlike of the undeformed case we are not able to obtain 
quantum deformation of Galilei group $G(1,3)$ by contraction of $SO_q(5)$.
It seems that quantum groups and corresponding quantum spaces are not a
suitable objects for simulation of noncommuting  space-time
because of the  fundamental physical correspondence principle is not satisfied
in this case.

\newpage
 \appendix
 \section{ Antipode $ S(T)=CT^{t}C^{-1}$
of quantum group  $ SO_v(N,\sigma,j) $ }

$$
S(t_{n+1,n+1})=t_{n+1,n+1}, \quad
S(t_{kk})=t_{kk}, \quad
S(\tau_{kk})=-\tau_{kk},
$$
$$
S(t_{k'k})=t_{k'k}\cosh 2J\rho_kv+
i\tau_{k'k}(\sigma_k,\sigma_{k'})\sinh 2J\rho_kv,
$$
$$
S(\tau_{k'k})=\tau_{k'k}\cosh 2J\rho_kv-
it_{k'k}(\sigma_k,\sigma_{k'})^{-1}\sinh 2J\rho_kv,
$$
$$
S(t_{k,n+1})=t_{n+1,k}\cosh J\rho_kv+
i\tau_{n+1,k}\displaystyle{\frac{(\sigma_{k'},\sigma_{n+1})}{(\sigma_k,\sigma_{n+1})}}
\sinh J\rho_kv,
$$
$$
S(\tau_{k,n+1})=\tau_{n+1,k}\cosh J\rho_kv-
it_{n+1,k}\displaystyle{\frac{(\sigma_{k},\sigma_{n+1})}{(\sigma_{k'},\sigma_{n+1})}}
\sinh J\rho_kv,
$$
$$
S(t_{n+1,k})=t_{k,n+1}\cosh J\rho_kv+
i\tau_{k,n+1}\displaystyle{\frac{(\sigma_{k'},\sigma_{n+1})}{(\sigma_k,\sigma_{n+1})}}
\sinh J\rho_kv,
$$
$$
S(\tau_{n+1,k})=\tau_{k,n+1}\cosh J\rho_kv-
it_{k,n+1}\displaystyle{\frac{(\sigma_{k},\sigma_{n+1})}{(\sigma_{k'},\sigma_{n+1})}}
\sinh J\rho_kv,
$$
\begin{eqnarray*}
S(t_{kp}) &=&  t_{pk}\cosh J\rho_kv\cosh J\rho_pv-
\frac{(\sigma_{p'},\sigma_{k'})}{(\sigma_k,\sigma_p)}
t'_{pk}\sinh J\rho_kv\sinh J\rho_pv+ \\
&&+i\frac{(\sigma_p,\sigma_{k'})}{(\sigma_k,\sigma_p)}
\tau_{pk}\sinh J\rho_kv\cosh J\rho_pv+ \\
&&+i\frac{(\sigma_{p'},\sigma_k)}{(\sigma_k,\sigma_p)}
\tau'_{pk}\cosh J\rho_kv\sinh J\rho_pv,
\end{eqnarray*}
\begin{eqnarray*}
S(t'_{kp}) & = & t'_{pk}\cosh J\rho_kv\cosh J\rho_pv-
t_{pk}\frac{(\sigma_p,\sigma_k)}{(\sigma_{k'},\sigma_{p'})}
\sinh J\rho_kv\sinh J\rho_pv- \\
& & -i\tau_{pk}\frac{(\sigma_p,\sigma_{k'})}{(\sigma_{k'},\sigma_{p'})}
\cosh J\rho_kv\sinh J\rho_pv-\\
& &-i\tau'_{pk}\frac{(\sigma_{p'},\sigma_k)}{(\sigma_{k'},\sigma_{p'})}
\sinh J\rho_kv\cosh J\rho_pv,
\end{eqnarray*}
\begin{eqnarray*}
S(\tau_{kp}) & = & \tau'_{pk}\cosh J\rho_kv\cosh J\rho_pv
+\tau_{pk}\frac{(\sigma_p,\sigma_{k'})}{(\sigma_k,\sigma_{p'})}
\sinh J\rho_kv\sinh J\rho_pv- \\
& & -it_{pk}\frac{(\sigma_p,\sigma_k)}{(\sigma_k,\sigma_{p'})}
\cosh J\rho_kv\sinh J\rho_pv+\\
 & & +it'_{pk}\frac{(\sigma_{p'},\sigma_{k'})}{(\sigma_k,\sigma_{p'})}
\sinh J\rho_kv\cosh J\rho_pv,
\end{eqnarray*}
\begin{eqnarray*}
S(\tau'_{kp}) & = & \tau_{pk}\cosh J\rho_kv\cosh J\rho_pv+
\tau'_{pk}\frac{(\sigma_{p'},\sigma_k)}{(\sigma_{k'},\sigma_p)}
\sinh J\rho_kv\sinh J\rho_pv- \\
& & -it_{pk}\frac{(\sigma_p,\sigma_k)}{(\sigma_{k'},\sigma_p)}
\sinh J\rho_kv\cosh J\rho_pv+\\
& & +it'_{pk}\frac{(\sigma_{p'},\sigma_{k'})}{(\sigma_{k'},\sigma_p)}
\cosh J\rho_kv\sinh J\rho_pv.
\end{eqnarray*}

 \section{  Coproduct  $ \Delta T = T \dot{\otimes} T $
of quantum group $ SO_v(N;\sigma;j)$}

\begin{eqnarray*}
\Delta t_{n+1,n+1} & = & t_{n+1} \otimes t_{n+1}+
2\sum^n_{k=1} [ (\sigma_k,\sigma_{n+1})^2t_{n+1,k}\otimes t_{k,n+1}+ \cr
 & &+(\sigma_{n+1},\sigma_{k'})^2\tau_{n+1,k} 
\otimes \tau_{k,n+1} ] ,
\end{eqnarray*}
\begin{eqnarray*}
 \Delta t_{kk} & = &t_{kk}\otimes t_{kk}+t_{k'k}\otimes t_{k'k}+
 (\sigma_k,\sigma_{n+1})^2t_{k,n+1} \otimes t_{n+1,k}+\\
& & +(\sigma_{n+1},\sigma_{k'})^2\tau_{k,n+1} \otimes \tau_{n+1,k}+ 
(\sigma_k,\sigma_{k'})^2(\tau_{k'k} \otimes \tau_{k'k}-
\tau_{kk} \otimes \tau_{kk})+ \\
& &+2\sum^n_{s=1, s\neq k}[(\sigma_k,\sigma_s)^2t_{ks}\otimes t_{sk}+ 
(\sigma_{k'},\sigma_{s'})^2t'_{ks}\otimes t'_{sk}+\\
& &+(\sigma_{k'},\sigma_s)^2\tau'_{ks}\otimes \tau_{sk}+
(\sigma_k,\sigma_{s'})^2\tau_{ks}\otimes \tau'_{sk}],
\end{eqnarray*}
\begin{eqnarray*}
\Delta \tau_{kk} & = &
\tau_{kk} \otimes t_{kk}+t_{kk} \otimes \tau_{kk}+
t_{k'k} \otimes \tau_{k'k}-\tau_{k'k} \otimes t_{k'k}+ \\
& & +\frac{(\sigma_k,\sigma_{n+1})(\sigma_{n+1},\sigma_{k'})}
{(\sigma_k,\sigma_{k'})}
(t_{k,n+1} \otimes \tau_{n+1,k}-\tau_{k,n+1} \otimes t_{n+1,k})+ \\
& &+\frac{2}{(\sigma_k,\sigma_{k'})}
\sum^n_{s=1, s\neq k}[(\sigma_k,\sigma_s)(\sigma_s,\sigma_{k'})
(t_{ks}\otimes \tau_{sk}-\tau'_{ks}\otimes t_{sk})+ \\
& &+(\sigma_k,\sigma_{s'})(\sigma_{s'},\sigma_{k'})
(\tau_{ks}\otimes t'_{sk}-t'_{ks}\otimes \tau'_{sk})],
\end{eqnarray*}
\begin{eqnarray*}
 \Delta t_{k'k}& = & t_{k'k} \otimes t_{kk}+t_{kk} \otimes t_{k'k}+
(\sigma_k,\sigma_{n+1})^2t_{k,n+1} \otimes t_{n+1,k}-\\
&&-(\sigma_{n+1},\sigma_{k'})^2\tau_{k,n+1} \otimes \tau_{n+1,k}+
(\sigma_k,\sigma_{k'})^2(\tau_{kk} \otimes \tau_{k'k}-\tau_{k'k} \otimes \tau_{kk})+\\
& &+2\sum^n_{s=1, s\neq k}[(\sigma_k,\sigma_s)^2t_{ks}\otimes t_{sk}- 
(\sigma_{k'},\sigma_{s'})^2t'_{ks}\otimes t'_{sk}+\\
& &+(\sigma_k,\sigma_{s'})^2\tau_{ks}\otimes \tau'_{sk}-
(\sigma_{k'},\sigma_s)^2\tau'_{ks}\otimes \tau_{sk}],
\end{eqnarray*}
\begin{eqnarray*}
\Delta\tau_{k'k} & = & \tau_{k'k} \otimes t_{kk}+t_{kk} \otimes \tau_{k'k}+
t_{k'k} \otimes \tau_{kk}-\tau_{kk} \otimes t_{k'k}+ \\
& & + \frac{(\sigma_k,\sigma_{n+1})(\sigma_{n+1},\sigma_{k'})}
{(\sigma_k,\sigma_{k'})}
(\tau_{k,n+1} \otimes t_{n+1,k}+t_{k,n+1} \otimes \tau_{n+1,k})+\\
& & +\frac{2}{(\sigma_k,\sigma_{k'})}
\sum^n_{s=1, s\neq k}[(\sigma_k,\sigma_s)(\sigma_s,\sigma_{k'})
(t_{ks}\otimes \tau_{sk}+
\tau'_{ks}\otimes t_{sk})+  \\
& & +(\sigma_k,\sigma_{s'})(\sigma_{s'},\sigma_{k'})(\tau_{ks}\otimes t'_{sk}
+t'_{ks}\otimes \tau'_{sk})],
\end{eqnarray*}
\begin{eqnarray*}
\Delta t_{k,n+1}& = & t_{k,n+1} \otimes t_{n+1,n+1}+
(t_{kk}+t_{k'k}) \otimes t_{k,n+1}+ \\
& &+\frac{(\sigma_k,\sigma_{k'})(\sigma_{k'},\sigma_{n+1})}
{(\sigma_k,\sigma_{n+1})}(\tau_{kk}+\tau_{k'k})\otimes \tau_{k,n+1}+ \\
& & +\frac{2}{(\sigma_k,\sigma_{n+1})}
\sum^n_{s=1, s \neq k}\biggl [ (\sigma_k,\sigma_s)(\sigma_s,\sigma_{n+1})
t_{ks}\otimes t_{s,n+1}+\\
& & + (\sigma_k,\sigma_{s'})(\sigma_{s'},\sigma_{n+1})\tau_{ks} \otimes \tau_{s,n+1}\biggr ],
\end{eqnarray*}
\begin{eqnarray*}
\Delta \tau_{k,n+1}& = & \tau_{k,n+1} \otimes t_{n+1,n+1}+
(t_{kk}-t_{k'k}) \otimes \tau_{k,n+1}+\\
& & +\frac{(\sigma_k,\sigma_{k'})(\sigma_k,\sigma_{n+1})}
{(\sigma_{n+1},\sigma_{k'})}(\tau_{k'k}-\tau_{kk})\otimes \tau_{k,n+1}+ \\
& & +\frac{2}{(\sigma_{n+1},\sigma_{k'})}
\sum^n_{s=1, s \neq k}\biggl [ (\sigma_{k'},\sigma_s)(\sigma_s,\sigma_{n+1})
\tau'_{ks}\otimes t_{s,n+1}+\\
& &+(\sigma_{k'},\sigma_{s'})(\sigma_{s'},\sigma_{n+1})t'_{ks} \otimes \tau_{s,n+1}\biggr ],
\end{eqnarray*}
\begin{eqnarray*}
\Delta t_{n+1,k}& = & t_{n+1,n+1} \otimes t_{n+1,k}+
t_{n+1,k}\otimes (t_{kk}+t_{k'k})+\\
& &+\frac{(\sigma_k,\sigma_{k'})(\sigma_{k'},\sigma_{n+1})}
{(\sigma_k,\sigma_{n+1})}\tau_{n+1,k}\otimes(\tau_{k'k}-\tau_{kk})+ \\
& & +\frac{2}{(\sigma_{k},\sigma_{n+1})}
\sum^n_{s=1, s \neq k}\biggl [ (\sigma_{k},\sigma_s)(\sigma_s,\sigma_{n+1})
t_{n+1,s}\otimes t_{sk}+\\
& & +(\sigma_{k},\sigma_{s'})(\sigma_{s'},\sigma_{n+1})\tau_{n+1,s}
\otimes \tau'_{sk}\biggr ],
\end{eqnarray*}
\begin{eqnarray*}
\Delta \tau_{n+1,k}& = & t_{n+1,n+1} \otimes \tau_{n+1,k}+
\tau_{n+1,k}\otimes (t_{kk}-t_{k'k})+\\
& &+\frac{(\sigma_k,\sigma_{k'})(\sigma_{k},\sigma_{n+1})}
{(\sigma_{n+1},\sigma_{k'})}t_{n+1,k}\otimes(\tau_{kk}+\tau_{k'k})+ \\
& & +\frac{2}{(\sigma_{n+1},\sigma_{k'})}
\sum^n_{s=1, s \neq k}\biggl [ (\sigma_{k'},\sigma_s)(\sigma_s,\sigma_{n+1})
t_{n+1,s}\otimes \tau_{sk}+\\
& &+(\sigma_{k'},\sigma_{s'})(\sigma_{s'},\sigma_{n+1})\tau_{n+1,s}
\otimes t'_{sk}\biggr ],
\end{eqnarray*}
\begin{eqnarray*}
\Delta t_{kp}& = & t_{kp} \otimes (t_{pp}+t_{p'p})
+(t_{kk}+t_{k'k})\otimes t_{kp}+\\
& &+\frac{(\sigma_k,\sigma_{n+1})(\sigma_{p},\sigma_{n+1})}
{(\sigma_k,\sigma_{p})}t_{k,n+1}\otimes t_{n+1,p}+ \\
& & +\frac{(\sigma_k,\sigma_{k'})(\sigma_{p},\sigma_{k'})}
{(\sigma_k,\sigma_{p})}(\tau_{kk}+\tau_{k'k})\otimes \tau'_{kp}+\\
& &+\frac{(\sigma_k,\sigma_{p'})(\sigma_{p},\sigma_{p'})}
{(\sigma_k,\sigma_{p})}\tau_{kp}\otimes (\tau_{p'p}-\tau_{pp})+ \\
& & +\frac{2}{(\sigma_{k},\sigma_{p})}
\sum^n_{s=1, s \neq k,p}\biggl [ (\sigma_{k},\sigma_s)(\sigma_s,\sigma_{p})
t_{ks}\otimes t_{sp}+\cr & & +
(\sigma_{k},\sigma_{s'})(\sigma_{s'},\sigma_{p})\tau_{ks}
 \otimes \tau'_{sp}\biggr ],
\end{eqnarray*}
\begin{eqnarray*}
\Delta t'_{kp}& = & t'_{kp} \otimes (t_{pp}-t_{p'p})
+(t_{kk}-t_{k'k})\otimes t'_{kp}+\\
& &+\frac{(\sigma_{n+1},\sigma_{k'})(\sigma_{n+1},\sigma_{p'})}
{(\sigma_{k'},\sigma_{p'})}\tau_{k,n+1}\otimes \tau_{n+1,p}+ \\
& & +\frac{(\sigma_k,\sigma_{k'})(\sigma_{k},\sigma_{p'})}
{(\sigma_{k'},\sigma_{p'})}(\tau_{k'k}-\tau_{kk})\otimes \tau_{kp}+\cr
& & + \frac{(\sigma_p,\sigma_{k'})(\sigma_{p},\sigma_{p'})}
{(\sigma_{k'},\sigma_{p'})}\tau'_{kp}\otimes (\tau_{pp}+\tau_{p'p})+ \\
& & +\frac{2}{(\sigma_{k'},\sigma_{p'})}
\sum^n_{s=1, s \neq k,p}\biggl [ (\sigma_{k'},\sigma_{s'})(\sigma_{s'},\sigma_{p'})
t'_{ks}\otimes t'_{sp}+ \cr & & +
(\sigma_{k'},\sigma_{s})(\sigma_{s},\sigma_{p'})\tau'_{ks}
 \otimes \tau_{sp}\biggr ],
\end{eqnarray*}
\begin{eqnarray*}
\Delta \tau_{kp}& = & \tau_{kp} \otimes (t_{pp}-t_{p'p})
+(t_{kk}+t_{k'k})\otimes \tau_{kp}+\\
& &+\frac{(\sigma_{p},\sigma_{n+1})(\sigma_{n+1},\sigma_{k'})}
{(\sigma_{k},\sigma_{p'})}\tau_{k,n+1}\otimes t_{n+1,p}+ \\
& & +\frac{(\sigma_k,\sigma_{k'})(\sigma_{k'},\sigma_{p'})}
{(\sigma_{k},\sigma_{p'})}(\tau_{kk}+\tau_{k'k})\otimes t'_{kp}+\cr & &+
\frac{(\sigma_k,\sigma_{p})(\sigma_{p},\sigma_{p'})}
{(\sigma_{k},\sigma_{p'})}t_{kp}\otimes (\tau_{pp}+\tau_{p'p})+ \\
& & +\frac{2}{(\sigma_{k},\sigma_{p'})}
\sum^n_{s=1, s \neq k,p}\biggl [ (\sigma_{k},\sigma_{s'})(\sigma_{s'},\sigma_{p'})
\tau_{ks}\otimes t'_{sp}+ \cr & & +
(\sigma_{k},\sigma_{s})(\sigma_{s},\sigma_{p'})t_{ks}
 \otimes \tau_{sp}\biggr ],
\end{eqnarray*}
\begin{eqnarray*}
\Delta \tau'_{kp} & = & \tau'_{kp} \otimes (t_{pp}+t_{p'p})
+(t_{kk}-t_{k'k})\otimes \tau'_{kp}+\\
 & &+\frac{(\sigma_{k},\sigma_{n+1})(\sigma_{n+1},\sigma_{p'})}
{(\sigma_{k'},\sigma_{p})}t_{k,n+1}\otimes \tau_{n+1,p}+ \\
 & & +\frac{(\sigma_k,\sigma_{k'})(\sigma_{k},\sigma_{p})}
{(\sigma_{k'},\sigma_{p})}(\tau_{k'k}-\tau_{kk})\otimes t_{kp}+ \cr & & +
\frac{(\sigma_{k'},\sigma_{p'})(\sigma_{p},\sigma_{p'})}
{(\sigma_{k'},\sigma_{p})}t'_{kp}\otimes (\tau_{p'p}-\tau_{pp})+ \\
& & +\frac{2}{(\sigma_{k'},\sigma_{p})}
\sum^n_{s=1, s \neq k,p}\biggl [ (\sigma_{k'},\sigma_{s})
(\sigma_{s},\sigma_{p})
\tau'_{ks}\otimes t_{sp}+\cr & &+
(\sigma_{k'},\sigma_{s'})(\sigma_{s'},\sigma_{p})t'_{ks}
 \otimes \tau'_{sp}\biggr ].
\end{eqnarray*}

 \section{ $(q-j)$-orthogonality relations $ TCT^t =C $
for quantum group $SO_v(N;\sigma;j) $}

Let us introduce the notation
$v_k=\rho_{k}v.$
\begin{eqnarray*}
1 & = & t^2_{n+1,n+1}+ 2\sum_{p=1}^{n} \Bigl\{
\left( (\sigma_p,\sigma_{n+1})^2t^2_{n+1,p} +
(\sigma_{n+1},\sigma_{p'})^2\tau^2_{n+1,p} \right)  \cosh{Jv_p}   +
\nonumber \\
& & +i(\sigma_p,\sigma_{n+1})(\sigma_{n+1},\sigma_{p'})
[t_{n+1,p},\tau_{n+1,p}]\sinh{Jv_p} \Bigr\},
\end{eqnarray*}
\begin{eqnarray*}
\frac{1}{2}\cosh{Jv_k} & = & (\sigma_k,\sigma_{n+1})^2t^2_{k,n+1}+
2\sum_{p=1,p \neq k}^{n} \Bigl\{
\left( (\sigma_k,\sigma_p)^2t^2_{kp}\right. + \\
&& + \left. (\sigma_k,\sigma_{p'})^2\tau^2_{kp} \right)  \cosh{Jv_p}   +
i(\sigma_k,\sigma_p)(\sigma_k,\sigma_{p'})
[t_{kp},\tau_{kp}]\sinh{Jv_p} \Bigr\} + \nonumber \\
& & +\frac{1}{2}\Bigl\{ t^2_{kk}+t^2_{k'k}+[t_{kk},t_{k'k}]_+ +
(\sigma_k,\sigma_{k'})^2
(\tau^2_{kk}+\tau^2_{k'k}+ \\
&&+[\tau_{kk},\tau_{k'k}]_{+}) \Bigr\}
\cosh{Jv_k} +\frac{i}{2}(\sigma_k,\sigma_{k'})
\Bigl\{ [t_{kk},\tau_{kk}]+ \\
&&+ [t_{k'k},\tau_{k'k}]+
[t_{k'k},\tau_{kk}]+[t_{kk},\tau_{k'k}] \Bigr\}\sinh{Jv_k},
\end{eqnarray*}
\begin{eqnarray*}
\frac{1}{2}\cosh{Jv_k} & = & (\sigma_{n+1},\sigma_{k'})^2\tau^2_{k,n+1}+
2\sum_{p=1,p \neq k}^{n} \Bigl\{
\left( (\sigma_{k'},\sigma_{p'})^2{t'}^2_{kp}\right. +\\
 && +\left.(\sigma_{k'},\sigma_{p})^2{\tau'}^2_{kp} \right)  \cosh{Jv_p}   -
i(\sigma_{k'},\sigma_{p'})(\sigma_{k'},\sigma_{p})
[t'_{kp},\tau'_{kp}]\sinh{Jv_p} \Bigr\}+  \nonumber \\
& & +\frac{1}{2} \Bigl\{ t^2_{kk}+t^2_{k'k}-[t_{kk},t_{k'k}]_{+}+
(\sigma_k,\sigma_{k'})^2(\tau^2_{kk}+\tau^2_{k'k}-\\
&&-[\tau_{kk},\tau_{k'k}]_{+}) \Bigr\}\cosh{Jv_k}
 + \frac{i}{2} (\sigma_k,\sigma_{k'})
\Bigl\{ [t_{kk},\tau_{kk}]+ \\
&&+[t_{k'k},\tau_{k'k}]-
[t_{k'k},\tau_{kk}]-[t_{kk},\tau_{k'k}] \Bigr\}\sinh{Jv_k},
\end{eqnarray*}
\begin{eqnarray*}
\frac{1}{2}\sinh{Jv_k} & = &
-i(\sigma_{k},\sigma_{n+1})(\sigma_{n+1},\sigma_{k'})t_{k,n+1}\tau_{k,n+1}+ \nonumber \\
& & +2\sum_{p=1,p \neq k}^{n} \Bigl\{
-i\left( (\sigma_{k},\sigma_{p})(\sigma_{k'},\sigma_{p})t_{kp}\tau'_{kp} + \right.\\
&&+\left.(\sigma_{k},\sigma_{p'})(\sigma_{k'},\sigma_{p'})\tau_{kp}t'_{kp} \right)
\cosh{Jv_p}+\left( (\sigma_{k},\sigma_{p})(\sigma_{k'},\sigma_{p'})
t_{kp}t'_{kp}\right.-\\
&&-\left.(\sigma_k,\sigma_{p'})(\sigma_{k'},\sigma_p)\tau_{kp}\tau'_{kp} \right) \sinh{Jv_p} \Bigr\}
 +\frac{i}{2}(\sigma_k,\sigma_{k'})
 \Bigl\{ [t_{kk},\tau_{kk}]-\\
&&-[t_{k'k},\tau_{k'k}]+[t_{k'k},\tau_{kk}]_{+}-
[t_{kk},\tau_{k'k}]_{+} \Bigr\}\cosh{Jv_k} + 
\frac{1}{2}\Bigl\{ t^2_{kk}-\\
&&-t^2_{k'k}+[t_{k'k},t_{kk}]+
(\sigma_k,\sigma_{k'})^2
\left( \tau^2_{kk}-\tau^2_{k'k}+[\tau_{k'k},\tau_{kk}]\right)
 \Bigr\}\sinh{Jv_k},
\end{eqnarray*}
\begin{eqnarray*}
\frac{1}{2}\sinh{Jv_k} & = &
i(\sigma_{k},\sigma_{n+1})(\sigma_{n+1},\sigma_{k'})\tau_{k,n+1}t_{k,n+1}+
 \nonumber \\
& & +2\sum_{p=1,p \neq k}^{n} \Bigl\{
i\left( (\sigma_{k'},\sigma_{p'})(\sigma_{k},\sigma_{p'})t'_{kp}\tau_{kp} 
\right.+\\
&&+\left.(\sigma_{k'},\sigma_{p})(\sigma_{k
},\sigma_{p})\tau'_{kp}t_{kp} \right)\cosh{Jv_p}+ 
\left( (\sigma_{k'},\sigma_{p'})(\sigma_{k},\sigma_{p})
t'_{kp}t_{kp}\right.-\\
&&-\left.(\sigma_k,\sigma_{p'})(\sigma_{k'},\sigma_p)\tau'_{kp}\tau_{kp}
 \right) \sinh{Jv_p} \Bigr\}+ 
\frac{i}{2}(\sigma_k,\sigma_{k'})
 \Bigl\{ [t_{kk},\tau_{kk}]-\\
&&-\left.[t_{k'k},\tau_{k'k}]-[t_{k'k},\tau_{kk}]_{+}+
[t_{kk},\tau_{k'k}]_{+} \Bigr\}\cosh{Jv_k}
 + \frac{1}{2}\Bigl\{ t^2_{kk}\right.-\\
&&-t^2_{k'k}-[t_{k'k},t_{kk}]
+(\sigma_k,\sigma_{k'})^2
\left( \tau^2_{kk}-\tau^2_{k'k}-[\tau_{k'k},\tau_{kk}]\right)
 \Bigr\}\sinh{Jv_k},
\end{eqnarray*}
\begin{eqnarray*}
0 & = &
(\sigma_{i},\sigma_{n+1})(\sigma_{k},\sigma_{n+1})t_{i,n+1}t_{k,n+1}+ \nonumber \\
& & +2\sum_{p=1,p \neq i,k}^{n} \Bigl\{
\left( (\sigma_{i},\sigma_{p})(\sigma_{k},\sigma_{p})t_{ip}t_{kp} +
(\sigma_{i},\sigma_{p'})(\sigma_{k},
\sigma_{p'})\tau_{ip}\tau_{kp} \right) \cosh{Jv_p}+ \nonumber \\
  & & +i\left( (\sigma_{i},\sigma_{p})(\sigma_{k},\sigma_{p'})
t_{ip}\tau_{kp}-(\sigma_k,\sigma_{p})(\sigma_{i},\sigma_{p'})\tau_{ip}t_{kp} \right) \sinh{Jv_p} \Bigr\}+  \nonumber \\
& & + \Bigl\{(\sigma_k,\sigma_{i})
  (t_{ii}+t_{i'i})t_{ki}+
  (\sigma_k,\sigma_{i'})(\sigma_i,\sigma_{i'})(\tau_{ii}+\tau_{i'i})
\tau_{ki} \Bigr\}\cosh{Jv_i} + \\
  & & +i \Bigl\{(\sigma_k,\sigma_{i'})
  (t_{ii}+t_{i'i})\tau_{ki}-
  (\sigma_i,\sigma_{i'})(\sigma_k,\sigma_{i})(\tau_{ii}+\tau_{i'i})t_{ki}
 \Bigr\}\sinh{Jv_i} + \\
& & +\Bigl\{ (\sigma_i,\sigma_k)t_{ik}(t_{kk}+t_{k'k})+
(\sigma_k,\sigma_{k'})(\sigma_i,\sigma_{k'})\tau_{ik}(\tau_{kk}+\tau_{k'k}) 
\Bigr\}\cosh{Jv_k} + \\
& & +i\Bigl\{-(\sigma_i,\sigma_{k'})\tau_{ik}(t_{kk}+t_{k'k})+
(\sigma_k,\sigma_{k'})(\sigma_i,\sigma_{k})t_{ik}(\tau_{kk}+\tau_{k'k})
 \Bigr\}\sinh{Jv_k},
\end{eqnarray*}
\begin{eqnarray*}
0 & = &
(\sigma_{i'},\sigma_{n+1})(\sigma_{k'},\sigma_{n+1})\tau_{i,n+1}\tau_{k,n+1}+ \nonumber \\
& & +2\sum_{p=1,p \neq i,k}^{n} \Bigl\{
\left( (\sigma_{i'},\sigma_{p'})(\sigma_{k'},\sigma_{p'})t'_{ip}t'_{kp} +
(\sigma_{i'},\sigma_{p})(\sigma_{k'},
\sigma_{p})\tau'_{ip}\tau'_{kp} \right) \cosh{Jv_p}+ \nonumber \\
  & & +i\left( -(\sigma_{i'},\sigma_{p'})(\sigma_{k'},\sigma_{p})
t'_{ip}\tau'_{kp}+(\sigma_{k'},\sigma_{p'})(\sigma_{i'},\sigma_{p})\tau'_{ip}t'_{kp} \right) \sinh{Jv_p} \Bigr\}+  \nonumber \\
& & + \Bigl\{ (\sigma_{k'},\sigma_{i'})(t_{ii}-t_{i'i})t'_{ki}+
  (\sigma_i,\sigma_{i'})(\sigma_{k'},\sigma_{i})(\tau_{i'i}-
\tau_{ii})\tau'_{ki} \Bigr\}\cosh{Jv_i} + \\
  & & +i\Bigl\{(\sigma_{k'},\sigma_{i})
  (t_{i'i}-t_{ii})\tau'_{ki}+
  (\sigma_i,\sigma_{i'})(\sigma_{k'},\sigma_{i'})(\tau_{i'i}-
\tau_{ii})t'_{ki} \Bigr\}\sinh{Jv_i} + \\
& & +\Bigl\{ (\sigma_{i'},\sigma_{k'})t'_{ik}(t_{kk}-t_{k'k})+
(\sigma_k,\sigma_{k'})(\sigma_{i'},\sigma_{k})\tau'_{ik}(\tau_{k'k}-
\tau_{kk}) \Bigr\}\cosh{Jv_k} + \\
& & +i \Bigl\{ (\sigma_{i'},\sigma_{k})\tau'_{ik}(t_{kk}-t_{k'k})+
(\sigma_k,\sigma_{k'})(\sigma_{i'},\sigma_{k'})t'_{ik}(\tau_{kk}-
\tau_{k'k}) \Bigr\}\sinh{Jv_k},
\end{eqnarray*}
\begin{eqnarray*}
0 & = &
-i(\sigma_{i},\sigma_{n+1})(\sigma_{k'},\sigma_{n+1})t_{i,n+1}\tau_{k,n+1}+ \nonumber \\
& & +2\sum_{p=1,p \neq i,k}^{n} \Bigl\{
-i\left( (\sigma_{i},\sigma_{p})(\sigma_{k'},\sigma_{p})t_{ip}\tau'_{kp} +
(\sigma_{i},\sigma_{p'})(\sigma_{k'},
\sigma_{p'})\tau_{ip}t'_{kp} \right) \cosh{Jv_p}+ \nonumber \\
  & & +\left( (\sigma_{i},\sigma_{p})(\sigma_{k'},\sigma_{p'})
t_{ip}t'_{kp}-(\sigma_{i},\sigma_{p'})(\sigma_{k'},\sigma_{p})\tau_{ip}\tau'_{kp} \right) \sinh{Jv_p} \Bigr\}+  \nonumber \\
& & +\Bigl\{ (\sigma_{k'},\sigma_{i'})
  (t_{ii}+t_{i'i})t'_{ki}-
  (\sigma_i,\sigma_{i'})(\sigma_{k'},\sigma_{i})(\tau_{ii}+
\tau_{i'i})\tau'_{ki} \Bigr\}\sinh{Jv_i} - \\
  & & -i\Bigl\{(\sigma_{k'},\sigma_{i})(t_{ii}+t_{i'i})\tau'_{ki}+
  (\sigma_i,\sigma_{i'})(\sigma_{k'},\sigma_{i'})(\tau_{ii}+\tau_{i'i})t'_{ki} 
\Bigr\}\cosh{Jv_i} + \\
& & +\Bigl\{ (\sigma_{i},\sigma_{k})t_{ik}(t_{kk}-t_{k'k})+
(\sigma_k,\sigma_{k'})(\sigma_{i},\sigma_{k'})\tau_{ik}(\tau_{kk}-
\tau_{k'k}) \Bigr\}\sinh{Jv_k} + \\
& & +i \Bigl\{(\sigma_{i},\sigma_{k'})\tau_{ik}(t_{k'k}-t_{kk})+
(\sigma_k,\sigma_{k'})(\sigma_{i},\sigma_{k})t_{ik}
(\tau_{kk}-\tau_{k'k}) \Bigr\}\cosh{Jv_k},
\end{eqnarray*}
\begin{eqnarray*}
0 & = &
i(\sigma_{k},\sigma_{n+1})(\sigma_{i'},\sigma_{n+1})\tau_{i,n+1}t_{k,n+1}+ \nonumber \\
& & +2\sum_{p=1,p \neq i,k}^{n} \Bigl\{
i\left( (\sigma_{i'},\sigma_{p'})(\sigma_{k},\sigma_{p'})t'_{ip}\tau_{kp} +
(\sigma_{k},\sigma_{p})(\sigma_{i'},
\sigma_{p})\tau'_{ip}t_{kp} \right) \cosh{Jv_p}+ \nonumber \\
  & & +\left( (\sigma_{i'},\sigma_{p'})(\sigma_{k},\sigma_{p})
t'_{ip}t_{kp}-(\sigma_{i'},\sigma_{p})(\sigma_{k},\sigma_{p'})\tau'_{ip}\tau_{kp} \right) \sinh{Jv_p} \Bigr\}+  \nonumber \\
& & + \Bigl\{ (\sigma_{k},\sigma_{i})(t_{ii}-t_{i'i})t_{ki}+
  (\sigma_i,\sigma_{i'})(\sigma_{k},\sigma_{i'})(\tau_{ii}-\tau_{i'i})\tau_{ki} 
\Bigr\}\sinh{Jv_i} + \\
  & & +i\Bigl\{(\sigma_{k},\sigma_{i'})(t_{ii}-t_{i'i})\tau_{ki}+
  (\sigma_i,\sigma_{i'})(\sigma_{k},\sigma_{i})(\tau_{i'i}-\tau_{ii})t_{ki} 
\Bigr\}\cosh{Jv_i}+ \\
& & + \Bigl\{ (\sigma_{i'},\sigma_{k'})t'_{ik}(t_{kk}+t_{k'k})-
(\sigma_k,\sigma_{k'})(\sigma_{i'},\sigma_{k})\tau'_{ik}(\tau_{kk}+\tau_{k'k})
 \Bigr\}\sinh{Jv_k} + \\
& & +i\Bigl\{(\sigma_{i'},\sigma_{k})\tau'_{ik}(t_{kk}+t_{k'k})+
(\sigma_k,\sigma_{k'})(\sigma_{i'},\sigma_{k'})t'_{ik}
(\tau_{kk}+\tau_{k'k}) \Bigr\}\cosh{Jv_k},
\end{eqnarray*}
\begin{eqnarray*}
0 & = &
(\sigma_{k},\sigma_{n+1})t_{k,n+1}t_{n+1,n+1}+
2\sum_{p=1,p \neq k}^{n} \Bigl\{
\left( (\sigma_{k},\sigma_{p})(\sigma_{p},\sigma_{n+1})t_{kp}t_{n+1,p}\right. +\\
&&+\left.(\sigma_{k},\sigma_{p'})(\sigma_{p'},
\sigma_{n+1})\tau_{kp}\tau_{n+1,p} \right) \cosh{Jv_p}+
i\left( (\sigma_{k},\sigma_{p})(\sigma_{p'},\sigma_{n+1})
t_{kp}\tau_{n+1,p} \right. -\\
&&-\left. (\sigma_{k},\sigma_{p'})(\sigma_{p},\sigma_{n+1})\tau_{kp}t_{n+1,p}
\right) \sinh{Jv_p} \left. \Bigr\}+\Bigl\{ (\sigma_{k},\sigma_{n+1})
  (t_{kk}+t_{k'k})t_{n+1,k}\right. + \\
&&+\left. (\sigma_k,\sigma_{k'})(\sigma_{k'},\sigma_{n+1})(\tau_{kk}+
\tau_{k'k})\tau_{n+1,k} \Bigr\}\cosh{Jv_k}\right. +  \\
& &+ i\Bigl\{(\sigma_{k'},\sigma_{n+1}) (t_{kk}+t_{k'k})\tau_{n+1,k}  - \\
& &-  (\sigma_k,\sigma_{k'})(\sigma_{k},\sigma_{n+1})(\tau_{kk}+
\tau_{k'k})t_{n+1,k} \Bigr\}\sinh{Jv_k},
\end{eqnarray*}
\begin{eqnarray*}
0 & = &
-i(\sigma_{k'},\sigma_{n+1})\tau_{k,n+1}t_{n+1,n+1}+ 
2\sum_{p=1,p \neq k}^{n} \Bigl\{
-i\left( (\sigma_{k'},\sigma_{p})(\sigma_{p},\sigma_{n+1})
\tau'_{kp}t_{n+1,p} \right. +\\
&&+\left.(\sigma_{k'},\sigma_{p'})(\sigma_{p'},
\sigma_{n+1})t'_{kp}\tau_{n+1,p} \right) \cosh{Jv_p}+
\left( -(\sigma_{k'},\sigma_{p'})(\sigma_{p},\sigma_{n+1})
t'_{kp}t_{n+1,p}\right. +\\
&& +\left. (\sigma_{k'},\sigma_{p})(\sigma_{p'},    
\sigma_{n+1})\tau'_{kp}\tau_{n+1,p} \right)  \sinh{Jv_p} \Bigr\}+ 
 \Bigl\{ (\sigma_{k},\sigma_{n+1}) (t_{k'k}-t_{kk})t_{n+1,k}+\\
&&+  (\sigma_k,\sigma_{k'})(\sigma_{k'},\sigma_{n+1})(\tau_{k'k}-
\tau_{kk})\tau_{n+1,k} \Bigr\}\sinh{Jv_k} + \\
& & +i\Bigl\{(\sigma_{k'},\sigma_{n+1})
(t_{k'k}-t_{kk})\tau_{n+1,k}-\\
&&-(\sigma_k,\sigma_{k'})(\sigma_{k},\sigma_{n+1})(\tau_{k'k}-\tau_{kk})t_{n+1,k}
\Bigr\}\cosh{Jv_k},
\end{eqnarray*}
\begin{eqnarray*}
0 & = &
(\sigma_{k},\sigma_{n+1})t_{n+1,n+1}t_{k,n+1}+ 
2\sum_{p=1,p \neq k}^{n} \Bigl\{ \Bigr.
\left( (\sigma_{k},\sigma_{p})(\sigma_{p},
\sigma_{n+1})t_{n+1,p}t_{kp}\right. +\\
&&+\left. (\sigma_{k},\sigma_{p'})(\sigma_{p'},
\sigma_{n+1})\tau_{n+1,p}\tau_{kp} \right) \cosh{Jv_p}+ 
i\left( (\sigma_{k},\sigma_{p'})(\sigma_{p},\sigma_{n+1})
t_{n+1,p}\tau_{kp}\right. -\\
&&-\left. (\sigma_{k},\sigma_{p})(\sigma_{p'},
\sigma_{n+1})\tau_{n+1,p}t_{kp} \right) \sinh{Jv_p} \Bigr\}+  
 \Bigl\{ (\sigma_{k},\sigma_{n+1})
  t_{n+1,k}(t_{kk}+t_{k'k})+\\
&&+  (\sigma_k,\sigma_{k'})(\sigma_{k'},\sigma_{n+1})\tau_{n+1,k}
(\tau_{kk}+\tau_{k'k}) \Bigr\}\cosh{Jv_k}+\\
&& +i\Bigl\{(\sigma_{k},\sigma_{k'})(\sigma_k,\sigma_{n+1})
  t_{n+1,k}(\tau_{kk}+\tau_{k'k})-\\
&&- (\sigma_{k'},\sigma_{n+1})\tau_{n+1,k}(t_{kk}+t_{k'k}) 
\Bigr\}\sinh{Jv_k},
\end{eqnarray*}
\begin{eqnarray*}
0 & = &
-i(\sigma_{k'},\sigma_{n+1})t_{n+1,n+1}\tau_{k,n+1}+
2\sum_{p=1,p \neq k}^{n} \Bigl\{
-i\left( (\sigma_{k'},\sigma_{p})(\sigma_{p},
\sigma_{n+1})t_{n+1,p}\tau'_{kp} \right. +\\
&&+\left. (\sigma_{k'},\sigma_{p'})(\sigma_{p'},
\sigma_{n+1})\tau_{n+1,p}t'_{kp} \right) \cosh{Jv_p}+ 
\left( (\sigma_{k'},\sigma_{p'})(\sigma_{p},\sigma_{n+1})
t_{n+1,p}t'_{kp}\right. -\\
&&-\left. (\sigma_{k'},\sigma_{p})(\sigma_{p'},\sigma_{n+1})\tau_{n+1,p}\tau'_{kp} 
\right) \sinh{Jv_p} \Bigr\}+  
 \Bigl\{ (\sigma_{k},\sigma_{n+1})t_{n+1,k}(t_{kk}-t_{k'k})+\\
&&+ (\sigma_k,\sigma_{k'})(\sigma_{k'},\sigma_{n+1})\tau_{n+1,k}
(\tau_{kk}-\tau_{k'k})\Bigr\}\sinh{Jv_k} +  \\
&&+i \Bigl\{(\sigma_{k},\sigma_{k'})(\sigma_k,\sigma_{n+1})
  t_{n+1,k}(\tau_{kk}-\tau_{k'k})-\\
&&-  (\sigma_{k'},\sigma_{n+1})\tau_{n+1,k}(t_{kk}-t_{k'k})
 \Bigr\} \cosh{Jv_k}.
\end{eqnarray*}

\end{document}